\documentclass[11pt,twoside]{article}
\usepackage{latexsym, amsfonts, amsmath}

\def\gR{\mathfrak g_{\mathbb R}}
\def\GR{G_{\mathbb R}}
\def\kR{\mathfrak k_{\mathbb R}}
\def\KR{K_{\mathbb R}}
\def\pR{\mathfrak p_{\mathbb R}}

\def\<{{\langle}}
\def\>{{\rangle}}

\def\ad{{\mathord{\rm ad}\,}}
\def\Ad{{\mathord{\rm Ad}}}
\def\Cent{ {\mathord{\rm Cent}}\,}

\def\Re{{\mathord{\rm Re}\,}}
\def\Im{{\mathord{\rm Im}\,}}
\def\Diff{\mathord{\rm  Diff}\thinspace}
\def\Lie{\mathord{Lie}}
\def\Vect{\mathfrak{Vect}\thinspace}
\def\Tr{\mathord{Tr}}
\def\Cinf{C^{\infty}}

\def\AND{\qquad\hbox{and}\qquad} 

\def\dimR{\dim_{\R}}
\def\dimC{\dim_{\C}}
\def\i{^{-1}}
\def\d{{\mathord{\rm d}}}
\def\({{\rm(}}  \def\){{\rm)}}
\def\half{\frac{1}{2}}

\def\rbot{\,\hbox to 5pt{\leaders\hrule\hfil}
   \vbox to 5pt{\leaders\vrule\vfil}\thinspace}

\def\pdb#1{\frac{\partial\phantom{#1}}{\partial#1}}
\def\ddt{\frac{\d}{\d t}}
\def\ddto{\frac{\d\phantom{t}}{\d{t}}\bigg|_{t=0}}

\def\del{\partial}
\def\odel{\ovl{\partial}}

\def\mapright#1{\thinspace\smash{\mathop{
\longrightarrow}\limits^{#1}}\thinspace}

\def\mapdown#1{\Big\downarrow\rlap{$\vcenter
{\hbox{$\scriptstyle#1$}}$}}

\def\C{{\mathbb C}}
\def\Z{{\mathbb Z}}
\def\R{{\mathbb R}}
\def\H{\mathbb H}
\def\RP{\R^{+}}

\def\cC{{\mathcal C}}
\def\cH{{\mathcal H}}
\def\cS{{\mathcal S}}
\def\cM{{\mathcal M}}
\def\cP{{\mathcal P}}

\def\A{{\mathcal A}}
\def\F{{\mathcal F}}
\def\N{{\mathcal N}}
\def\L{{\mathcal L}}
\def\T{{\mathcal T}}
\def\Q{{\mathcal Q}}
\def\V{{\mathcal V}}

\def\a{{\mathfrak a}}
\def\b{{\mathfrak b}}

\def\g{{\mathfrak g}}

\def\k{{\mathfrak k}}

\def\p{{\mathfrak p}}

\def\s{{\mathfrak s}}

\def\u{{\mathfrak u}}
\def\v{{\mathfrak v}}

\def\x{{\mathfrak x}}
\def\y{{\mathfrak y}}

\def\so{{\mathfrak s\mathfrak o}}

\def\sl{{\mathfrak s\mathfrak l}}

\def\su{\s\u}

\def\ovl{\overline}

\def\al{\alpha}
\def\be{\beta}
\def\ep{\varepsilon}  \def\tep{\tilde{\varepsilon}}
\def\ga{\gamma}
\def\ka{\kappa}
\def\la{\lambda} 
\def\om{\omega}  \def\Om{\Omega}
\def\sig{\sigma}  \def\Sig{\Sigma} 
\def\th{\theta} \def\vth{\vartheta} \def\Th{\Theta}
\def\ze{\zeta}

\def\oz{{\ovl{z}}}

\def\of{\ovl{f}}

\def\ov{{\ovl{v}}}

\def\I{{\mathbf I}}  \def\J{{\mathbf J}}  \def\K{{\mathbf K}}
\def\bi{{\mathbf i}} \def\bj{{\mathbf j}} \def\bk{{\mathbf k}}
\def\omI{\om_{\I}} \def\omJ{\om_{\J}}\def\omK{\om_{\K}}

\def\O{{\mathcal O}}  \def\OR{{\mathcal O}_{\R}}

\def\rhoo{\rho_o}
\def\cs{{\mathcal R}} 
\def\Uvth{U^{\vth}}
\def\Rhol{R^{\mathord{\rm hol}}}
\def\dc{\d^{c}}
\def\Him{\H^{im}}
\def\vs{\v^\sharp}
\def\vC{\v_{\C}}
\def\Mk{{\cM(\kappa)}}
\def\x{{\mathfrak a}}
\def\y{{\mathfrak b}}
 \def\vsig{\varsigma}

\def\tbe{\tilde{\be}}
\def\Zr{Z{(\R)}}  \def\Ur{U{(\R)}}
\def\qq{_{q}}

\def\phir{\phi{(\R)}}
\def\VR{V_{\mathbb R}}
\def\XR{X_{\mathbb R}}
 
\def\GRR{G(\R)}
\def\XRR{X(\R)}
\def\oZ{\ovl{Z}}

\newtheorem{theorem}{Theorem}[section]
\newtheorem{proposition}[theorem]{Proposition}
\newtheorem{lemma}[theorem]{Lemma}
\newtheorem{corollary}[theorem]{Corollary}
 
\newtheorem{example}[theorem]{Example}       
\newtheorem{remark}[theorem]{Remark}         

\newenvironment{proof}{{\noindent {\bf Proof}\,\,}}
{\hfill$\Box$\medskip} 

\title{\bfseries Instantons and  Kaehler Geometry
of   Nilpotent Orbits}
\author{Ranee BRYLINSKI 
\thanks{Research supported in part  by NSF  Grant No. DMS-9505055} 
\\[5pt]
\itshape Department of Mathematics\\
\itshape Pennsylvania State University\\
\itshape University Park, PA 16802\\
\itshape USA\\[12pt]}
\date{}

\begin{document}           
\maketitle

\begin{abstract}
The first obstacle  in building a Geometric Quantization theory
for    nilpotent orbits of a real
semisimple Lie group  has been the 
 lack of an  invariant polarization. In order to generalize the 
Fock  space construction of the quantum mechanical 
oscillator, a polarization 
of the symplectic orbit invariant  under the maximal compact
subgroup  is required.
 
In this paper, we explain how such a polarization on the  
orbit arises naturally  from the work of Kronheimer and
Vergne.  This occurs in the context of hyperkaehler geometry.
The polarization is complex and in fact makes the orbit into a 
(positive) Kaehler manifold.
We study  the geometry of this Kaehler structure, the Vergne
diffeomorphism,  and the Hamiltonian functions giving the
symmetry.   We indicate how all this fits into a   quantization
program.
\end{abstract}

\section{Introduction}\label{sec:1}
\setcounter{equation}{0}

Quantization is a procedure for constructing a quantum
system with symmetry out  of a classical system with
symmetry.  No axiomatic or even systematic method of
quantization is known.  Instead, quantization exists as an
empirical science, made up of a growing series of examples. 
In many ways, quantization is an art.

The nature of the classical and quantum systems under consideration depends on 
the context and on the scope of the investigation.  At present a universal sort of
quantization scheme seems completely out of reach.  Such a scheme would have to
include the quantization of gravity as well as the quantization of classical field
theory.  In fact, some physicists believe that even familiar classical
theories  must be modified in order that they can be  ``quantized" to
give a consistent and meaningful quantum theory. 
 
There is a rather clear ``beginning level" at which to
formulate and study the  quantization problem.   This is the
case where one starts with a classical Hamiltonian
mechanical (dynamical) system with  symmetry.  Such a
system is given by a phase space $(M,\om)$ together with a
Hamiltonian function $F$ and  a Lie subalgebra 
\begin{equation}\label{eq:11}\g \subset C^{\infty}(M)
\end{equation}   Here $(M,\om)$ is a real
symplectic manifold; i.e.,    $\om$ is a symplectic form on a 
smooth manifold $M$  of dimension $2n$.  The symplectic
structure $\om$ defines a Poisson bracket
$\{\,,\,\}$  on
$C^{\infty}(M)$, giving it the structure of a Poisson algebra.  
The Hamiltonian is a fixed smooth function $F:M \rightarrow
\R$ which determines the time evolution (the dynamics) of the
system.  Each smooth function $\phi$ on $M$ determines a
Hamiltonian vector field
$\xi_{\phi}$ by the formula $\xi_{\phi}\rbot\om+d\phi=0$. 
The Poisson  bracket is given by 
$\{ \phi,\psi \} = \xi_{\phi}(\psi) = \om(\xi_{\phi},\xi_{\psi})$
and satisfies
\begin{equation}\label{eq:12}
[\xi_{\phi},\xi_{\psi}] = \xi_{\{ \phi,\psi \}}\end{equation}

Thus the Hamiltonian vector fields $\xi_{\phi}$ of the functions 
$\phi\in\g$  give an infinitesimal Lie algebra action of $\g$ on
$M$.  This constitutes \emph{infinitesimal  Hamiltonian 
symmetry}.  If the $\g$-action integrates to a smooth action of
a  Lie  group $G$ on $M$, then this $G$-action is called
Hamiltonian.  Regardless of integration, the inclusion  
$\g\hookrightarrow C^{\infty}(M)$ defines  a smooth
infinitesimally
$\g$-equivariant  \emph{moment map}
\begin{equation}\label{eq:13}
 M\rightarrow \g^*\end{equation} 
 
If  $G$-acts transitively on $M$ in a Hamiltonian fashion, then
the moment map (\ref{eq:13}) is just a covering onto a
coadjoint orbit of $G$. This focuses attention on coadjoint
orbits as the ``elementary Hamiltonian spaces" (up to
covering).  Moreover, each coadjoint orbit  $P$  has a
canonical symplectic structure $\sig$, sometimes called the
KKS (Kirillov-Kostant-Souriau) symplectic structure, derived 
from the Lie algebra bracket on $\g$.  Indeed, each $x\in\g$
defines a linear function $\phi^x$ on $\g^*$ and hence on $P$.
Then $\sig$ is the 
unique symplectic form such that the mapping
\begin{equation}\label{eq:14} 
\g\to\Cinf(P),\qquad x\mapsto\phi^x
\end{equation}
preserves brackets, i.e., is a Lie algebra homomorphism.
This discussion applies equally well in the category of
holomorphic symplectic manifolds; cf. \cite{bkHam}.

An outstanding problem is the quantization of conical
coadjoint orbits $\OR$ of a real semisimple Lie group $\GR$.
Here we take $\GR$ to be a real form of a connected and
simply-connected complex semisimple Lie group $G$
with Lie algebra $\g$.  We
may identify $\gR=\Lie\GR$ with its dual
$\gR^*$, and then the conical orbits $\OR$ identify with the
orbits of nilpotent elements in $\gR$. These are the so-called
``nilpotent orbits". The irreducible unitary representations
arising from quantization of nilpotent orbits  are often called
``unipotent" representations. These are examples of
``singular representations".

In order to start a quantization program for nilpotent orbits, 
we need at the outset an invariant polarization of $\OR$.
In analogy with the well-known quantization of the  harmonic
oscillator, we want a polarization invariant under a maximal
compact subgroup $\KR$ of $\GR$. (In general $\OR$
does not admit $\GR$-invariant polarizations.)

Remarkably, a $\KR$-invariant polarization arises, in a 
uniform natural manner on every real nilpotent orbit, from
the work of Kronheimer and Vergne. This comes about by first
working on the complexification $\O$ of $\OR$; $\O\subset\g$
is a complex nilpotent orbit of  $G$. We let $\I$ denote the
natural complex structure on $\O$; then $\OR$ is an $\I$-real
form.

Kronheimer (\cite{Kr2})  in 1990 identified each complex
nilpotent orbit $O$ as an instanton moduli space. In
particular   the holomorphic symplectic structure $(\I,\Sig)$
on $O$ extends to a   hyperkaehler structure
$(g,\I,\J,\K,\omI,\omJ,\omK)$ where $\Sig=\omJ+i\omK$. 
We outline the Kronheimer model of $\O$    in \S\ref{sec:5}.

Then  Vergne in 1995 used this to discover a diffeomorphism
\begin{equation}\label{eq:15} 
\V:\OR\to Y
\end{equation}
of  each real  nilpotent orbit $\OR$ to a complex  
$K$-homogeneous cone $Y$, where $K\subset G$ is the
complexification of $\KR$.  This recovered the
Kostant-Sekiguchi (\cite{sek}) correspondence. 
 
The upshot of  Vergne's work on the Kronheimer 
instanton model of $\O$ is that $\OR$ is a $\J$-complex
submanifold of $\OR$. Moreover, $\OR$ is then a Kaehler
submanifold of $\O$ with respect to $(\J,\omJ)$. But then
$\omJ=\Re\Sig$ is just the real KKS symplectic form on $\OR$.
 
So the ``new" complex structure $\J$ provides a complex 
polarization on $\OR$ and moreover makes $\OR$ into a
Kaehler manifold which identifies with $Y$ as a complex
manifold.  We explain the Vergne theory and the
Kaehler structure in \S\ref{sec:6} and  \S\ref{sec:7}.
In \S\ref{sec:6} we also give a different proof of Vergne's
result (see especially Proposition \ref{prop:67} and
Corollaries \ref{cor:66} and \ref{cor:69}).
 
In \S\ref{sec:7},\ref{sec:8},\ref{sec:9} we develop the 
properties of the Vergne diffeomorphism and the Kaehler
structure.  Our first main result is the Triple Sum formula in 
in Theorem \ref{thm:78}. This leads to our key result for
quantization in Theorem \ref{thm:94} on the nature of the
Hamiltonian functions $\phi^x$. Indeed we have the Cartan
decomposition $\gR=\kR\oplus\pR$. While the Hamiltonian
flows of the functions $\phi^x$, $x\in\kR$, preserve $\J$, the
flows of  the remaining functions $\phi^v$, $v\in\pR$, do not 
preserve $\J$. The question then is how will they quantize. On
a  classical level, we can ask how to write down $\phi^v$ in
terms of  holomorphic and antiholomorphic  functions. The
answer in Theorem \ref{thm:94} is that $\phi^v$ is the real
part of a holomorphic function. The interpretation  of this is
discussed further in \S\ref{sec:9}.

In \S\ref{sec:8}, we explain another aspect of the Kaehler
structure, namely that there is a global Kaehler potential  
$\rhoo$ on $(\OR,\J,\sig)$. This function $\rhoo:\OR\to\R$
is $\KR$-invariant and uniquely determined by the condition
that it is homogeneous of degree $1$  under the Euler scaling
action of $\RP$. 

This Kaehler potential arises by restriction from the
hyperkaehler potential on $\O$. In
\S\ref{sec:2},\S\ref{sec:3},\S\ref{sec:4} we develop the  basic
theory of hyperkaehler manifolds,
hyperkaehler cones and hyperkaehler potentials
based on  results from \cite{HKLR} and \cite{Sw}.

The importance of the Kaehler potential $\rhoo$ is this: in
our quantization program for real nilpotent orbits, $\rhoo$
plays the role of the Hamiltonian, i.e.,  the energy function.
Moreover, $\rhoo$ gives rise in  Theorem    \ref{thm:87} 
and Corollary \ref{cor:87}  to a
realization of $T^*Y$ as a holomorphic symplectic
complexification of $\OR$.

Our quantization program building on this geometry
will be developed in subsequent papers.
See also  \cite{B1},\cite{B2}, \cite{bkLagr}. 
 
In the quantization of $(\OR,\J,\sig)$, we want to ``quantize"
the Hamiltonian  functions $\phi^z$, $z\in\gR$,  by converting
the $\phi^z$ into self-adjoint operators $\Q(\phi^z)$ on a
space of holomorphic half-forms on $(\OR,\J)\cong  Y$. 
The conversion must satisfy in particular Dirac's axiom that
Poisson bracket of functions goes over into the commutator of
operators so that 
$\Q(\{\phi^z,\phi^w\})=i[\Q(\phi^z),\Q(\phi^w)]$.

A main idea coming out of Corollary \ref{cor:87} is that we can
try  to ``promote" the functions  $\phi^z$ on $\OR$ to rational 
functions on the holomorphic symplectic complexification
$T^*Y$. For this to work, we need some sort of
analyticity and algebraicity for  the embedding of $\OR$ into 
$T^*Y$. 

The appropriate  notion combining analyticity and algebraicity
here turns out to be that of  a Nash embedding.  In the
Appendix, we give an outline of Nash geometry, starting from
the theory of real algebraic varieties.  O. Biquard has proven in
[Bi] that the hyperkaehler potential on $\O$, and hence the
$SO(3)$-action on $\O$ and Vergne diffeomorphism
(\ref{eq:15}), are Nash.

I thank Alex Astashkevich, Olivier Biquard, Nestor Handzy, Bert Kostant,
Michele Vergne, and Francois Ziegler for useful conversations.
I also thank Nestor Handzy for help in writing this paper.
Parts of this work were carried out during visits to other institutions and I thank
them for their hospitality:  Institute for Advanced Study (Spring 1995 term),
Harvard University (summers of 1995 and 1996), and Brown University (summer  
of 1997).

\section{Hyperkaehler Manifolds}\label{sec:2}
\setcounter{equation}{0}

In this section, we review and perhaps clarify some basic
notions of hyperkaehler geometry that we use
throughout this paper.

A hyperkaehler manifold
$(X,g,\J_1,\J_2,\J_3)$ is real manifold $X$ 
of  dimension $4n$ together with a Riemannian metric
$g$ and three  complex
structures $\J_1$,$\J_2$,$\J_3$ such that
(i) $\J_1\J_2\J_3=-1$ and (ii) $g$ is  a Kaehler metric with respect
to each of $\J_1$,$\J_2$,$\J_3$.

Then by (i), $\J_1$,$\J_2$,$\J_3$
satisfy the quaternion relations 
\[\J_1^2=\J_2^2=\J_3^2=-1,\qquad \J_a\J_b=\ep_{abc}\J_c\]
Here $a,b,c\in\{1,2,3\}$ are distinct and
$\ep_{abc}=\hbox{sgn}(abc)$. Thus every 
tangent space of $X$ becomes a   quaternionic vector space.  

By (ii), $X$ has three Kaehler manifold structures 
$(\J_1,\om_1)$,  
$(\J_2,\om_2)$, $(\J_3,\om_3)$, all with Kaehler metric $g$.
The  Kaehler forms $\om_a$ are given  by 
$g(u,v)=\om_a(u,\J_av)$. 
We call these Kaehler manifolds  $X_1$, $X_2$, $X_3$,
respectively. 
 
The data $(X,g,\om_1,\om_2,\om_3)$
serves equally well to define the hyperkaehler structure as  we
may recover the complex structures by the formula
\begin{equation}\label{eq:221} 
\om_c(u,v)=\om_a(\J_bu,v)\ep_{abc}
\end{equation}

We define   three complex $2$-forms on $X$
\begin{equation}\label{eq:222}
\Om_1=\om_2+i\om_3,\hskip 10pt \Om_2=\om_3+i\om_1, 
\hskip 10pt\Om_3=\om_1+i\om_2\end{equation}
Then $\Om_a$ is $\J_a$-holomorphic. This is shown in     
\cite[pp. 549-550]{HKLR}.

Inside the quaternion algebra 
\[\H=\R\oplus\R\bi\oplus\R\bj\oplus\R\bk\]
we have the standard $2$-sphere
\[S^2=\{q=a\bi+b\bj+c\bk\,|\, |q|=1\}\]
of pure imaginary quaternions of unit norm.

Corresponding  to   a point $q=a\bi+b\bj+c\bk$ on  $S^2$, we
have the pair
\[ \J\qq=a\J_1+b\J_2+c\J_3 \AND
\om\qq=a\om_1+b\om_2+c\om_3\]
Then $(X,g,\J_q,\om_q)$ is again a Kaehler structure on $X$
with complex structure $\J\qq$ and Kaehler form $\om\qq$; 
we write $X\qq$ for this  Kaehler manifold.
Thus we have a $2$-sphere $\cS_X$ of Kaehler structures  
$(\J\qq,\om\qq)$   on $X$ and we have identified  
$\cS_X$  with $S^2$.

Let $q\mapsto \tau q$ be the standard rotation action of 
$SO(3)$  on $S^2$. This induces an  $SO(3)$-action on $\cS_X$
given by 
$\tau\cdot\J\qq=\J_{\tau q}$ and 
$\tau\cdot\om\qq=\om_{\tau q}$. Let
\begin{equation}\label{eq:223} 
\cC\qq\subset SO(3) 
\end{equation}
be the circle subgroup of  which fixes  $q\in S^2$.  

The  generalization of (\ref{eq:222}) is  that any $q'\in S^2$ 
orthogonal to $q$ determines a $\J\qq$-holomorphic
symplectic form    
$\om_{q'}+i\om_{q''}$  on $X$
where   $q''=q\times q'$  is   the cross product of $q$ with $q'$.
 
\begin{example}\label{ex:23} \rm
The first  example of a   hyperkaehler manifold  is  the 
flat  quaternionic vector space. Let $X=\R^{4n}$ with standard
linear coordinates $x^r_s$ where $r=0,1,2,3$ and  $s=1,\dots,n$.
We may make $X$ into an $n$-dimensional quaternionic vector
space,  where   $\H$ acts by left multiplication,
in the obvious way so that the functions 
\begin{equation}\label{eq:231} 
q_s=x^0_s+x^1_s\bi+x^2_s\bj+x^3_s\bk
\end{equation}
are quaternionic linear coordinates. 

The following data defines a hyperkaehler structure on $X$:  left
multiplication by $\bi$, $\bj$ and $\bk$ give  the complex 
structures $\J_1,\J_2,\J_3$ so that
\begin{equation}\label{eq:232} 
\J_a\left(\pdb{x^0_s}\right)=\pdb{x^a_s},\qquad  
\J_a\left(\pdb{x^b_s}\right)=\pdb{x^c_s}
\end{equation} 
where $(abc)$ is a cyclic permutation of $1,2,3$.
Under $\J_a$, $X$ identifies with $\C^{2n}$ with   linear 
holomorphic coordinate functions $x^0_s+ix^a_s$, 
$x^b_s+ix^c_s$.   The hyperkaehler metric is
\begin{equation}\label{eq:233}
g=\sum_{r,s}(\d x^r_s)^2
\end{equation}
The   three Kaehler  forms $\om_1,\om_2,\om_3$ are
\begin{equation}\label{eq:234}
\om_a=\sum_{s=1}^n\d x^0_s\wedge\d x^a_s+
\d x^b_s\wedge\d x^c_s\end{equation}
\end{example}

Next we introduce hyperkaehler symmetry into the 
picture.  Let $U$ be a Lie group. A \emph{hyperkaehler action} of  
$U$  on $(X,g,\J_q,\om_q)$ is a smooth Lie group action  of $U$
on $X$   which preserves all the hyperkaehler structure.

From now on we assume that  $U$ is a compact connected 
semisimple Lie group and we have a   hyperkaehler action of $U$
on $X$.  Then differentiation gives an infinitesimal   action 
of  the Lie algebra $\u$ of $U$ by the vector fields $\xi^u$  where
$\xi^u_p=\frac{\d}{\d t}\big|_{t=0}(\exp -tu)\cdot p$ at
$p\in X$. In other words,   we get a Lie algebra
homomorphism
\begin{equation}\label{eq:241} 
\u\to\Vect X, \qquad u\mapsto\xi^u
\end{equation}

Now consider each Kaehler  manifold $X\qq$. We let
$\Cinf(X)_{\om\qq}$  denote the algebra $\Cinf(X)$ equipped
with the Poisson bracket defined by $\om\qq$.

The $U$-action  on $X$ is  symplectic with respect to
$\om\qq$ and consequently, since $\u$ is semisimple, is
Hamiltonian.    This means that we can   solve the
equations $\xi^u\rbot\om_q+\d\ze\qq^u=0$ for functions 
$\ze\qq^u$ such that $\{\ze\qq^u,\ze\qq^v\}=\ze\qq^{[u,v]}$.
The  momentum functions $\ze\qq^u$ are uniquely determined.
So we get a Lie algebra homomorphism
\begin{equation}\label{eq:242}
\u\to \Cinf(X)_{\om\qq},\hskip 1.5pc u\mapsto\ze\qq^u
\end{equation} 
The corresponding  $U$-invariant  moment map 
\begin{equation}\label{eq:243} 
\ze\qq:X\to\u\end{equation}
is defined by  $\ze\qq^u(p)=(u,\ze\qq(p))_\u$.
Here we  identify $\u\simeq\u^*$ by  means of the Killing form 
$(~,~)_{\u}$.

Consider now the three moment maps 
$\ze_1=\ze_{\bi}$, $\ze_2=\ze_{\bj}$, $\ze_3=\ze_{\bk}$,
Putting  these together we obtain a triple moment map
\begin{equation}\label{eq:244}
\ze=(\bi\ze_1,\bj\ze_2,\bk\ze_3):
X\to \bi\u\oplus \bj\u\oplus \bk\u
\end{equation}

Let $G$ be the complexification of $U$. 
Then $G$ is the complex semisimple algebraic group 
characterized   by either of the following properties:
(i) any linear representation of $U$ on a complex 
(finite-dimensional) vector space extends uniquely 
to a  linear  representation of $G$, or
(ii) $U$ is a compact real form of $G$.  It follows from (ii)
that $U$ and $G$ have the same fundamental group.

We  assume now that $U$, and hence $G$, is simply-connected.
The Lie algebra  of $G$ is the  complex semisimple Lie algebra
\begin{equation}\label{eq:251}
\g=\u\otimes\C=\u\oplus \bi\u
\end{equation}
We identify $\g\simeq\g^*$ using   the complex Killing form 
$(~,~)_{\g}$ on $\g$.  
We note that $(u,v)_{\u}=(u,v)_{\g}$ for $u,v\in\u$.
This follows because   an $\R$-linear  map $L:\u\to\u$
determines  a $\C$-linear map $L_\C:\g\to\g$ and then 
$\Tr_\R L=\Tr_\C L_\C$. 
 
Now we consider the holomorphic symplectic manifolds
$(X,\J_a,\Om_a)$, $a=1,2,3$. We let $\Rhol(X_a)$ denote the
algebra of
$\J_a$-holomorphic functions on $X$ equipped with the 
Poisson bracket defined by $\Om_a$. 

Since $\L_{\xi^u}\J_a=\J_a$, $u\in\u$, it follows that
$\half\xi^u$ is  the real part of a unique $\J_a$-holomorphic
vector field  $\Xi_a^u$ on $X$. Precisely, 
$\half\xi^u=\Re\Xi^u_a$ where
\begin{equation}\label{eq:252}
\Xi_a^u=\half(\xi^u-i\J_a\xi^u)
\end{equation}
Then we get the bracket relations for $u,v\in\u$
\[\{\Xi_a^u,\Xi_a^v\}=\Xi_a^{[u,v]}\]

Now we have an  infinitesimal $\J_a$-holomorphic Lie algebra action 
\begin{equation}\label{eq:253}
\g\to\Vect_{\J_a-hol}(X),\hskip 2pc
z=u+iv\mapsto\Xi_a^z=\Xi_a^u+i\Xi_a^v
\end{equation}
of $\g$ on  $X_a$. This is the complexification of the 
infinitesimal $\u$-action   (\ref{eq:241}).   If (\ref{eq:252}) integrates to a
holomorphic $G$-action on $(X,\J_a)$, then we will
say that the $U$-action on $X$ \emph{complexifies} with respect
to 
$\J_a$. 

Regardless of   integration, the infinitesimal 
action (\ref{eq:253}) of $\g$ preserves the 
holomorphic symplectic form $\Om_a$ defined in (\ref{eq:222}).
Then, since $\g$ is semisimple,  the infinitesimal action (\ref{eq:253}) is 
Hamiltonian.  We
have  a unique complex Lie algebra homomorphism
\begin{equation}\label{eq:254}
\g\to \Rhol(X_a),\hskip 1.5pc z\mapsto\Phi_a^z
\end{equation}
given by momentum functions $\Phi_a^z$ so that
$\Xi_a^z\rbot\Om_a+\d\Phi_a^z=0$.
Then
\begin{equation}\label{eq:255}
\Phi_a^{u+iv}=\Phi_a^u+i\Phi_a^v\AND
\Phi^{u}_a=\ze^u_{b}+i\ze^u_{c}
\end{equation}
where $(abc)$ is a cyclic permutation of $1,2,3$.

The corresponding 
$\g$-equivariant  $\J_a$-holomorphic moment map
\begin{equation}\label{eq:256}
\Phi_a:X_a\to\g
\end{equation}
is defined by  $\Phi^{z}_a(p)=(z,\Phi_a(p))_\g$.
Thus we get  the three maps 
\begin{equation}\label{eq:257}\Phi_1=\ze_2+i\ze_3,\hskip 1.5pc
\Phi_2=\ze_3+i\ze_1,\hskip 1.5pc
\Phi_3=\ze_1+i\ze_2\end{equation}

The formulas (\ref{eq:255})-(\ref{eq:257}) encode a lot of information about
the coupling of the complex and symplectic structures on $X$,
as $\Phi_a$ is $\J_a$-holomorphic.
In particular they show how the real functions 
$\ze_1^u$, $\ze_2^u$, $\ze_3^u$
give rise to holomorphic functions on $X$.

\begin{example}\label{ex:26} \rm
We continue  the discussion of $X=\H^n$ from Example \ref{ex:23}.
Let $U$ be the group of all $\R$-linear transformations of $X$ 
which  preserve $g$ and commute with the  $\H^*$-action on $X$. Then  $U$
is the  familiar model of  the compact symplectic group $Sp(n)$.
Clearly  this $U$-action preserves all the hyperkaehler data on
$\H^n$.  In the case $n=1$ then $U\simeq SU(2)$ and moreover  
$U$  acts by right multiplication by quaternions of  unit norm.

The $U$-action complexifies, with respect to any 
complex structure  $\J\qq\in\cS_{X}$, to a complex linear
complex algebraic action of $G\simeq Sp(2n,\C)$ on $\H^n$.
This action is transitive on $\H^n-\{0\}$.

We have a free $\Z_2$-action  on $\H^n-\{0\}$ by
multiplication by $\pm 1$. This $\Z_2$-action
preserves all the hyperkaehler data on $\H^n-\{0\}$
and commutes with the $SU(2)$ and $G$-actions.
The quotient $\O=(\H^n-\{0\})/\Z_2$ inherits  a
hyperkaehler structure, an action of  
$SU(2)/\Z_2\simeq SO(3)$ and a $U$-action.
\end{example}

\section{ Hyperkaehler  Cones}\label{sec:3}
\setcounter{equation}{0}

In this section we  explain the notion of a hyperkaehler cone.

To begin with we recall that a \emph{symplectic cone} of 
positive integer weight $k$    is a symplectic manifold
$(M,\om)$  together with a  smooth action  
\[\RP\times M\to M,\qquad (t,m)\mapsto\ga_t(m)\]  
of the group $\RP$ of positive real numbers such that  
\begin{equation}\label{eq:321}
\ga_{t}^*\om=t^k\om\end{equation}
This means that the $\RP$-action scales the symplectic form 
and it  has  weight $k$. 

The prototype example is the case where  $k=1$ and  
 $M=T^*Q$   is a  cotangent bundle with its canonical symplectic
structure. Here  $\RP$ acts on $T^*Q$  by the linear scaling 
action on the fibers of the projection $T^*Q\to Q$.

Let $\eta$ be the  infinitesimal generator  of the $\RP$-action.
Then differentiating (\ref{eq:321}) we get the equivalent
condition
\begin{equation}\label{eq:322}
\L_\eta\om=k\om\end{equation}  
It follows  that $\om$ is exact  with symplectic potential 
$\frac{1}{k}(\eta\rbot\om)$;   i.e., 
\[\om=\d(\frac{1}{ k}\eta\rbot\om)\] 
We conclude in particular that a  symplectic cone  is 
non-compact  (and has positive dimension).

  Next we define a \emph{Kaehler cone} of weight $k$ to be a 
Kaehler manifold $(Z,\J,\om,g)$ together with a  smooth action 
\begin{equation}\label{eq:331}
\ga:\C^*\times Z\to Z, \qquad
(s,m)\mapsto\ga_s(m)\end{equation} 
which satisfies the  three conditions
\begin{equation}\label{eq:332} 
(i) \mbox{ the action $\ga$ is holomorphic}, \qquad  
(ii) \ga_s^*\om=|s|^k\om,\qquad    (iii) \ga_s^* g=|s|^kg
\end{equation} 
(These  are consistent  with  redundancy).

The condition (i) means that the map (\ref{eq:331}) is 
holomorphic. So (i) implies  $\ga_s^*\J=\J$. Also any two of 
$\ga_s^*\J=\J$, (ii), (iii) imply the other.

To work out (ii) and (iii), we  use the product decomposition
\[\C^*=\RP\times S^1\]
So the $\C^*$ action splits into a product of an $\RP$-action with
an $S^1$-action.  Then (ii) and (iii) say:  $\om$ and $g$ are 
homogeneous of  degree $k$ under the  $\RP$-action, but they
are fixed by the $S^1$-action.  

Thus  $S^1$ acts by Kaehler
automorphisms.   In particular, the $S^1$-action is 
symplectic and so  has a  moment map on $Z$ at least locally
with values in $\R$. We can write this moment map as
$\frac{k}{2}\rho$.  Then
$\rho$ is  a local Kaehler potential, i.e., $\rho$ satisfies
\begin{equation}\label{eq:333}  
i\del\odel\rho=\om\end{equation}
where $d=\del+\odel$ is the standard
decomposition of $\d$ into $(1,0)$ and $(0,1)$ parts.

The Lie algebra of $\C^*$ is $\C=\R\oplus\R\bi$ with 
$[1,\bi]=0$. Differentiating the $\C^*$-action we get an
infinitesimal vector  field action
\begin{equation}\label{eq:334}
\psi:\C\to\Vect Z,\qquad v\mapsto 
\psi^v=\ddto\ga_{\exp -tv}\end{equation} 
We put
\begin{equation}\label{eq:335}
\eta=\psi^{-1} \AND \th=\psi^{-\bi}
\end{equation}
so that $\eta$ and $\th$ are, respectively, the infinitesimal
generators of  the actions of $\RP$ and $S^1$.  Since $\psi$ 
is a  (real) Lie algebra homomorphism we have
\begin{equation}\label{eq:336}[\eta,\th]=0\end{equation}

Notice that the infinitesimal generator of the holomorphic
action of $\C^*$ on $Z$ is the holomorphic vector field
\begin{equation}\label{eq:337}E=\half\eta-\half i\th\end{equation}
Now, we can give an equivalent infinitesimal version of the  
conditions (\ref{eq:332}) on  $\ga$:
$\psi$ must be  $\C$-linear, i.e.,
\begin{equation}\label{eq:338}\th=\J\eta\end{equation}
and also
\begin{equation}\label{eq:339}
\L_\eta\J=0,\quad\L_\eta\om=k\om,\quad\L_\eta g=k g,
\quad\L_{\th}\,\J=\L_{\th}\,\om=\L_{\th}\,g=0\end{equation}
Notice that the condition $\L_\eta\J=0$ itself implies  
$[\eta,\J\eta]=0$.

Now we define a \emph{hyperkaehler cone} of weight $k$ to be a  
hyperkaehler manifold
$(X,g,\J_q,\om_q)$ together with a left  $\H^*$-action 
\begin{equation}\label{eq:341}\ga:\H^*\times X\to X, \qquad
(h,m)\mapsto\ga_h(m)\end{equation}
which satisfies the three conditions
\begin{equation}\label{eq:342} 
\begin{array}{ccl}(i)& \ga_h^*\J_q &=\;\J_{h\i qh}
\mbox{ and the action of $\C^*_q$ on $X$ is 
$\J_q$-holomorphic}  \\
(ii)& \ga_h^*\om_q &=\; |h|^k\om_{h\i qh}\\
(iii)& \ga_h^* g\phantom{_q}&=\; |h|^kg 
\end{array}  
\end{equation}
Here  $h\in\H^*$,  $q\in S^2$ and
\begin{equation}\label{eq:343}
\C^*_q=\{a+bq\,|\, (a,b)\in\R^2, (a,b)\neq(0,0)\}
\end{equation}
Again (i)-(iii)  are consistent  with redundancies.  
It suffices to check (i)-(iii) just for $q=\bi,\bj,\bk$.

We have the natural direct product decomposition
\[\H^*=\RP\times SU(2)\]
The formulation of  (i)-(iii)  in terms of the component actions 
of  $\RP$ and $SU(2)$ is:  
(a) for all $q$, $\RP$   acts $\J_q$-holomorphically and
scales  $g$ and $\om_q$ so that they have weight $k$, 
(b) the action of  the circle 
\begin{equation}\label{eq:344}
T^1_q=\{\cos t+q\sin t\,|\, t\in\R\}\subset SU(2)\end{equation}
on $X$ is $\J_q$-holomorphic, and (c)
$SU(2)$ acts isometrically and  permutes the   Kaehler 
structures  $X_q$ according to the   standard action  of 
$SU(2)=\widetilde{SO}(3)$ on $S^2$.

We can also rewrite the conditions (i)-(iii) at the infinitesimal 
level. Indeed the Lie algebra of the multiplicative group 
$\H^*$ is
$\H$ with Lie bracket $[u,v]=uv-vu$.  We have a standard basis
$\bj_0,\bj_1,\bj_2, \bj_3$ of $\H$ with 
\[\bj_0=1,\quad\bj_1=\bi,\quad\bj_2=\bj,\quad\bj_3=\bk\]
and bracket relations
\begin{equation}\label{eq:345}[\bj_0,\bj_a]=0\AND
[\bj_a,\bj_b]=2\ep_{abc}\bj_c\end{equation}

Differentiating the $\H^*$-action we get an infinitesimal 
vector  field action
\begin{equation}\label{eq:346}
\psi:\H\to\Vect X,\qquad q\mapsto\psi^q\end{equation}
We put
\begin{equation}\label{eq:347}
\eta=\psi^{-1}\AND \th_a=\psi^{-\bj_a},\quad a=1,2,3\end{equation}
Then $\eta$ is  the infinitesimal 
generator of the 
$\RP$-action and $\th_a$ is the infinitesimal generator of the 
$T^1_q$-action. The bracket relations are
\begin{equation}\label{eq:348}
[\eta,\th_a]=0,\qquad [\th_a,\th_b]=-2\ep_{abc}\th_c
\end{equation}

Now we can give an equivalent version of the conditions 
(\ref{eq:342}) on $\ga$:  $\psi$ is $\H$-linear, i.e.,
\begin{equation}\label{eq:349}
\th_a=\J_a\eta,\quad a=1,2,3
\end{equation}
and also
\begin{equation}\label{eq:3410}
\begin{array}{lll}
\phantom{|}\L_{\eta}\om_a =k\om_a, &
\phantom{|}\L_{\eta}\J_a=0,& \phantom{|}\L_{\eta}g=kg\\
\L_{\th_{a}}\om_{b}=-2\tep_{abc}\om_{c},\phantom{XX}&
\L_{\th_{a}}\J_{b}=-2\tep_{abc}\J_{c},\phantom{XX} &
\L_{\th_{a}}g=0
\end{array}
\end{equation}
In (\ref{eq:3410}), we do not require $a,b,c$ distinct, but 
instead  we define 
\[\tep_{abc}=\left\{  
\begin{array}{ll}
sgn(abc)& \mbox{if $a,b,c$ are distinct}\\ 
0& \mbox{otherwise}
\end{array}\right.     \]
We call a vector field  action (\ref{eq:346}) satisfying 
(\ref{eq:347})-(\ref{eq:3410})
\emph{infinitesimally conical} with weight  $k$.

\begin{lemma}\label{lem:34} 
Any infinitesimally  conical vector field  action  of $\H$ 
on $X$ necessarily has weight $k=2$.
\end{lemma}

\begin{proof}  If $a,b,c$ are distinct then (\ref{eq:349}) and 
(\ref{eq:221}) give
\begin{equation}\label{eq:3411}
\th_a\rbot\om_b=-\ep_{abc}\,\eta\rbot\om_c
\end{equation}
But then
\[-2\ep_{abc}\om_{c}=\L_{\th_{a}}\om_{b}=
\d(\th_a\rbot\om_b)=
-\ep_{abc}\d(\eta\rbot\om_c)=-\ep_{abc}\L_\eta\om_c=
-k\ep_{abc}\om_c\]  Hence $k=2$.
\end{proof}

So  any hyperkaehler cone necessarily has weight $2$.
From now on, we assume $k=2$ in (\ref{eq:342}) and 
(\ref{eq:3410}).

\section{The Hyperkaehler  Potential}\label{sec:4}
\setcounter{equation}{0}

In this section, we explain, on the global level, the
relation between the hyperkaehler cone structure and
the   hyperkaehler potential. This was worked out locally  in
\cite{Sw}; see also   \cite[pg. 553]{HKLR} for part of this.

Corresponding to each complex structure
$\J\qq$ on $X$ , we have the decomposition 
$\d=\del\qq+\odel\qq$
of the exterior derivative into $(1,0)$ and $(0,1)$ parts.
We put 
\begin{equation}\label{eq:421}
\dc\qq\,=\,-\half\J_q\d\,=\,-\frac{i}{ 2}
(\del\qq-\odel\qq)\end{equation}
A global \emph{hyperkaehler potential} on $X$ is a smooth  
function  $\rho: X\to\R$ which is a simultaneous 
Kaehler potential for each Kaehler structure  $X\qq$,  i.e.,
\begin{equation}\label{eq:422}
\om\qq=i\del\qq\odel\qq\rho=\d\dc\qq\rho
\end{equation}
for all $q\in S^2$.
It follows easily that $\rho$ is a hyperkaehler potential
iff $\rho$ is a Kaehler potential for $X_1$, $X_2$ and $X_3$.

\begin{proposition}\label{prop:42} 
Suppose $(X,\J_q,\om_q,g)$  admits a global 
hyperkaehler potential   $\rho:X\to\R$. 
Let $\eta$ be the  vector field  on $X$ defined by 
\begin{equation}\label{eq:423} 
\eta\rbot g=\d\rho
\end{equation}  and set $\th_q=\J_q\eta$. 
 
Then $\rho$, after perhaps being modified by  adding
a constant, satisfies $\eta\rho=2\rho$ so that $\rho$ is
homogeneous  of weight  $2$. 
The vector fields $\eta,\th_1,\th_2,\th_3$
define an infinitesimally conical vector field action of
$\H$ on $X$ where  $\th_1,\th_2,\th_3$ define the 
infinitesimal $\so(3)$-action.
 
The potential $\rho$ is $SO(3)$-invariant. For each  $q\in
S^2$,  the Hamiltonian flow of $\rho$ with respect to $\om_q$
integrates  the vector field $\th_q$. In other words, we have
\begin{equation}\label{eq:424}
\th_q\rbot\om_q+\d\rho=0\end{equation}
so that $\rho$ is a simultaneous moment map for each
infinitesimal $S^1$-action  defined by   $\th_q$. 
\end{proposition}
\begin{proof}
To begin with, we observe that if 
a vector field  $\eta$  and a function $\rho$ satisfy 
(\ref{eq:423}) then
\begin{equation}\label{eq:425}
\mbox{$\rho$ is a hyperkaehler potential} 
\quad\iff\quad
\mbox{$\L_{\eta}\om_a=2\om_a$ for $a=1,2,3$}\end{equation}
Indeed    the computation
\[\<\eta\rbot\om_a,\xi\>=\om_a(\eta,\xi)=g(\J_a\eta,\xi)
=-g(\eta,\J_a\xi)=\<-\d\rho,\J_a\xi\>=
\<-\J_a\d\rho,\xi\>=\<2\dc_a\rho,\xi\>\] 
gives $\eta\rbot\om_a=2\dc_a\rho$. Applying $\d$ to this we
get 
\[\L_{\eta}\om_a=2\d\dc_a\rho\]
This implies (\ref{eq:425}).
 
Now suppose $\rho$ is a hyperkaehler potential.
We need to prove all the  relations in (\ref{eq:3410}); in fact we can
ignore the two involving the metric $g$ as they are redundant.
The   relations $\L_{\eta}\om_a=2\om_a$  are done in  (\ref{eq:425}).
These imply the relations  $\L_{\eta}\J_b=0$ because we can 
apply $\L_\eta$ to   (\ref{eq:221}). This in turn gives 
$[\eta,\th_a]=\L_\eta(\J_a\eta)=0$. To prove the second line of 
relations in (\ref{eq:3410})  we first observe that
$g(u,v)=\om_a(u,\J_av)=\om_a(-\J_au,v)$ gives
\begin{equation}\label{eq:426}
-\th_a\rbot\om_a=\eta\rbot~g\end{equation}
So 
\[\L_{\th_a}\om_a=\d(\th_a\rbot\om_a)= 
-\d(\eta\rbot g)=\d^2\rho=0\]
For  $a,b,c$  distinct we  find using (\ref{eq:3411})  
\[\L_{\th_a}\om_b=\d(\th_a\rbot\om_b)=
-\ep_{abc}\d(\eta\rbot\om_c)=-\ep_{abc}\L_\eta\om_c=
-2\ep_{abc}\om_c\]
Next, applying   $\L_{\th_c}$ to  (\ref{eq:221}) we find that
$\L_{\th_c}\J_a=-2\ep_{abc}\J_b$.
This proves the six independent  relations in (\ref{eq:3410}). Also 
\[[\th_a,\th_b]=\L_{\th_a}\th_b=\L_{\th_a}(\J_b\eta)
=-2\ep_{abc}\J_c\eta=-2\ep_{abc}\th_c\]
Thus we have an infinitesimally conical $\H$-action.

Now   applying  $\L_\eta$ to (\ref{eq:423}) we find 
$\L_\eta(\d\rho)=2\d\rho$. So    $\eta\d\rho=2\d\rho$ and
therefore $\eta\rho=2\rho+C$ where $C$ is constant.  We
can replace $\rho$ by  $\rho+C/2$ so that $\eta\rho=2\rho$.
Next applying  $\L_{\th_a}$ to (\ref{eq:423}) we find
$\L_{\th_a}(\d\rho)=0$.  So 
$\L_{\th_a}\rho=C_a$ where $C_a$ is constant.
But then it follows from the semisimplicity of $\su(2)$, in 
particular  from the   relations
$[\th_a,\th_b]=-2\ep_{abc}\th_c$, that  $C_1=C_2=C_3=0$. Hence
$\th_a\rho=0$. 

Finally, it suffices to check (\ref{eq:424}) for any single  $q\in S^2$
because of the $\so(3)$-action.  Clearly
$\th_a\rbot\om_a+\d\rho=0$ follows by (\ref{eq:423}) and 
(\ref{eq:426}).
\end{proof}

We may think of a hyperkaehler $\H^*$-conical structure  
on $X$ as a family of Kaehler cone structures on $X$ which satisfy
additional properties. Indeed, for each $q\in S^2$, $X_q$ is a 
Kaehler cone with respect to the action of
$\C^*\qq=\RP\times\cC\qq$ where the action of $\cC\qq$
integrates $\th_q$.
 
A converse to Proposition \ref{prop:42} follows easily from the proof.

\begin{corollary}\label{cor:43} 
Suppose $H^1_{deRham}(X)=0$ and
$X$ admits an infinitesimally conical vector field action of $\H$. 
Then {\rm(\ref{eq:423})}
has a smooth solution $\rho$ on $X$ {\rm(}unique
up to the addition of a constant function{\rm)}, and  $\rho$ is
a global hyperkaehler potential.  The further
condition $\eta\rho=2\rho$ uniquely determines $\rho$.
\end{corollary}
\begin{proof}
It is enough because of  (\ref{eq:425}) and 
the relation $\L_\eta\om_a=2\om_a$ in (\ref{eq:3410}) to produce a
solution $\rho$ to  (\ref{eq:423}).
Since $H^1_{deRham}(X)=0$, the problem reduces to showing 
that the $1$-form $\eta\rbot g$ is closed. 
This is easy:  the general fact (\ref{eq:426})  and one of the relations in
(\ref{eq:3410})  give $\d(\eta\rbot~g)=-\L_{\th_a}\om_a=0$.  
\end{proof}

\begin{example}\label{ex:44} \rm
We continue discussing $X=\H^n$ from Examples \ref{ex:23} and
\ref{ex:26}. The left multiplication action of $\H^*$ on $X$ given by
\[h\circ(q_1,\dots,q_n) =(hq_1,\dots,hq_n)\]
makes $X$ into a
hyperkaehler cone of weight $2$. The natural $SO(3)$-action on 
the $2$-sphere  $\cS_{X}$ of Kaehler  structures is induced by  
the $SU(2)$-action on $X$ defined by
left multiplication by quaternions of  unit norm. 
The infinitesimal generator $\eta$ of the $\RP$-action and the 
corresponding hyperkaehler potential $\rho$ are
\begin{equation}\label{eq:441}\eta=\sum_{r,s}x^{r}_{s}
\frac{\partial\phantom{X}}{\partial x^{r}_{s}}
\AND
\rho=\half\sum_{r,s}(x^r_s)^2=
\half\sum_{s}|q_s|^2\end{equation}
The vector fields $\th_a$, $a=1,2,3$,  are
\begin{equation}\label{eq:442}
\th_a=x^0_s\pdb{x^a_s}-x^a_s\pdb{x^0_s}+
x^b_s\pdb{x^c_s}-x^c_s\pdb{x^b_s}
\end{equation}
where $(abc)$ is cyclic.
\end{example}

Let
$\Him=\bi\R\oplus\bj\R\oplus\bk\R$.   We have a weight-$2$ 
left action of $\H^*$ on $\Him$ defined by
\begin{equation}\label{eq:451}
h\bullet w=|h|^2(hwh\i)\end{equation}
This is just the product of the degree
$2$ scaling action of $\RP$ with the  spin $1$ action of  $SU(2)$. 
We then get a tensor product action of  $\H^*\times U$ on 
\begin{equation}\label{eq:452} 
\bi\u\oplus\bj\u\oplus\bk\u=\Him\otimes\u 
\end{equation}
given by, for $w\in\Him$ and $u\in\u$,
\begin{equation}\label{eq:453}
(h,a)\cdot(wu)=(h\bullet w)\Ad_au\end{equation}

\begin{lemma}\label{lem:45} 
Suppose  a hyperkaehler cone  $X$ has a hyperkaehler
action of $U$ which commutes with  the $\H^*$-action. Then
the triple moment map {\rm(\ref{eq:244})} is  equivariant under    
$\H^*\times U$ with respect to the  action {\rm(\ref{eq:453})}.
\end{lemma}
\begin{proof}
This follows immediately by transforming the equation
$\xi^u\rbot\om_a+\d\ze_{a}^u=0$ under the $\H^*$-action.
\end{proof}

\begin{lemma}\label{lem:46} 
Let $Z$ be a complex manifold with $H^{1}_{DeRham}Z=0$.  Then
\item{\rm (i)} A smooth function $\rho:Z\to\R$ is pluriharmonic 
{\rm(}i.e., $\del\odel\rho=0${\rm)} if and only if
$\rho$ is the real part of a holomorphic function on $Z$.
\noindent Suppose further a compact group
$H$ acts holomorphically on $Z$ and this
action complexifies to a transitive action of  the
complexified group $H_\C$ on $Z$.  Then 
\item{\rm (ii)}The only $H$-invariant 
pluriharmonic  smooth functions $\rho:X\to\R$ 
are the   constants.
\item{\rm (iii)}  If $\om_Z$ is a Kaehler form on $Z$,
then an  $H$-invariant Kaehler potential on $Z$,
if it exists, is unique up to addition of a constant.
\end{lemma}
\begin{proof}
(i) This follows by a  standard argument (see
e.g. \cite{Kra}).  Indeed the $1$-form $i(\del-\odel)(\rho)$
is  closed and so exact since $H^{1}_{DeRham}Z=0$. So 
we can solve $i(\del-\odel)(\rho)=\d\phi$
globally  for a smooth real-valued function $\phi$. Then
$i\del\rho=\del\phi$ and $i\odel\rho=i\odel\phi$.
But then $f=\rho-i\phi$ satisfies $\odel f=0$
which means $f$ is holomorphic.
(ii) By (i), we have  $\rho=\Re f$ for
some holomorphic function $f$. We see easily that $f$ is
$H$-invariant. But then $f$ is $H_\C$-invariant. Now the 
transitivity  assumption   forces $f$, and hence
$\rho$, to be constant. (iii) follows from (ii) since the difference
of any two potentials is pluriharmonic.
\end{proof}

An immediate consequence is 

\begin{proposition}\label{prop:46} 
Suppose $X$ is as in Proposition {\rm \ref{prop:42}} and $X$
carries a hyperkaehler action of $U$ which commutes with the
infinitesimal $\H$-action. Then the homogeneous degree $2$
solution $\rho$ to  {\rm(\ref{eq:423})} is, up to addition of a constant, 
the unique $U$-invariant hyperkaehler potential on $X$.  
\end{proposition}

\section{Instantons and the Kronheimer  Model 
of a Complex Nilpotent Orbit}\label{sec:5}
\setcounter{equation}{0}

In this section we  recall some main results from  \cite{Kr2}.
First  we construct  the instanton space $\Mk$  , and then we 
explain how $\Mk$ is isomorphic to a complex nilpotent orbit
$\O$.

The subspace  $\Him$  of pure imaginary quaternions
is a Lie subalgebra of $\H$ isomorphic to
$\so(3)$; cf. (\ref{eq:345}).  From now on, we may identify $\so(3)$ 
with $\Him$.  The adjoint action of $SO(3)$ identifies with
the  standard spin $1$ action of 
$SO(3)$ on $\Him$, which we write as 
$(\tau,w)\mapsto\tau*w$.
 
Let $L(\so(3),\u)$ be the space of  $\R$-linear  maps 
\begin{equation}\label{eq:521} 
A:\so(3)\to\u, \qquad w\mapsto A_w
\end{equation} 
An element $A(t)$ of the
path space $\cP=C^\infty(\R,L(\so(3),\u))$ is given by
a triple  
\[A(t)=(A_1(t),A_2(t),A_3(t))\] where  $A_a(t)=A_{-\bj_a}(t)$.
The adjoint actions of $SO(3)$ and $U$ define
a representation of $SO(3)\times U$ on $L(\so(3),\u)$. This gives 
a  representation of $SO(3)\times U$ on $\cP$ defined by
\begin{equation}\label{eq:522} 
(\tau,u)\cdot A_w(t)=\Ad_u\cdot
A_{\tau\i*w}(t) 
\end{equation}

Let $\cM\subset\cP$ be the subspace of 
paths $A(t)$ which satisfy
the system of  three differential equations
\begin{equation}\label{eq:523} 
\ddt A_a=-2A_a-[A_b,A_c] 
\end{equation}
where $(abc)$ is a cyclic permutation of $1,2,3$.
Now let  
\begin{equation}\label{eq:524} 
\ka:\so(3)\to\u 
\end{equation}
be a non-zero Lie algebra homomorphism; then  $\ka$ is
$1$-to-$1$.  We put
\[e_a=\ka(-\bj_a),\qquad a=1,2,3\]

Let $\Mk\subset\cM$ be the subspace
of paths satisfying the boundary conditions
\begin{equation}\label{eq:525} 
\lim_{t\to\infty} A(t)=0,\hskip 20pt
\lim_{t\to -\infty}A(t)\in C(\ka) 
\end{equation} 
where  $C(\ka)\subset L(\so(3),\u)$ is the  space
of Lie algebra homomorphisms which are
$U$-conjugate to $\ka$. The action of $SO(3)\times U$ on $\cP$
leaves stable $\cM$ and $\Mk$.  Thus we have a smooth action
\begin{equation}\label{eq:526} 
SO(3)\times U\to\Diff \Mk 
\end{equation}

An \emph{instanton} is an element  of $\Mk$. We often write an
instanton   $A(t)$ simply as ``$A$'' with the time dependence
being understood.

It follows from the equations (\ref{eq:523}) (see 
\cite{Kr2}) that we have three well-defined 
$U$-equivariant  mappings
\begin{equation}\label{eq:531} 
\ze_1,\ze_2,\ze_3:\Mk\to\u 
\end{equation}
given by
\begin{equation}\label{eq:532} 
\ze_a(A)=\frac{1}{2}\lim_{t\to\infty}e^{2t}A_a(t)
\end{equation}
 
We can  arrange these  three maps $\ze_1,\ze_2,\ze_3$ into the
single map
\begin{equation}\label{eq:533} 
\ze=(\bi\ze_1,\bj\ze_2,\bk\ze_3):\Mk\to
\bi\u\oplus\bj\u\oplus\bk\u 
\end{equation}
We have a natural (tensor product) action of $SO(3)\times U$ on 
$\bi\u\oplus\bj\u\oplus\bk\u$ because of (\ref{eq:452}). We see easily
that  $\ze$ is $(SO(3)\times U)$-equivariant.

Each triple $d=(d_1,d_2,d_3)$
lying in $C(\ka)$ gives rise to  the ``model"  instanton
\[D(t)=(1+e^{2t})\i(d_1,d_2,d_3)\]
Then $\ze_a(D(t))=\frac{1}{2}d_a$.

It is easy to check that we have an action of
$\RP$ on $\Mk$ given by
\begin{equation}\label{eq:541} 
\la\diamond A(t)=A\big(t-\half\log\la\big) 
\end{equation}
This commutes with the action of $SO(3)\times U$.

The complexification of $\ka$ is a complex Lie algebra 
embedding 
\[\ka_\C:\so(3,\C)\to\g\]
The set  $\N$ of  all nilpotent elements in $\so(3,\C)$ is
\[\N=\{a\bi+b\bj+c\bk\,|\,  a^2+b^2+c^2=0\}\]
and $\N-\{0\}$ is  a single orbit under
the adjoint action of $SO(3,\C)$. It follows that  all the 
elements
$\ka(z)$, $z\in\N-\{0\}$, lie in a single adjoint orbit $\O$
in $\g$. Then $\O$ consists of nilpotent elements in $\g$;
i.e.,  $\O$ is a complex nilpotent orbit.

In particular then
\begin{equation}\label{eq:551} 
\O=G\cdot(e_2+ie_3)=G\cdot(e_3+ie_1)=G\cdot(e_1+ie_2)
\end{equation}
Then  $\O$ inherits a complex   structure $\I$   from
the natural embedding of $\O$ into $\g$; we call this
embedding
\begin{equation}\label{eq:552} 
\Phi_\I:\O\to\g 
\end{equation}
Or, equivalently,  $\I$ is induced by the $G$-action.
$\I$ and the  $G$-invariant KKS  holomorphic symplectic 
form $\Sigma$ make $\O$ into a    holomorphic symplectic
manifold 
\begin{equation}\label{eq:553} 
(\O,\I,\Sigma) 
\end{equation}
 
We recall that $\Sigma$ is the unique holomorphic symplectic   
form   on $\O$ such that the adjoint action of $G$ on $\O$ is
holomorphic Hamiltonian with moment map $\Phi_\I$.
In terms of   the holomorphic component  function $\Phi_\I^z$,
$z\in\g$, defined by $\Phi_\I^z(w)=(z,w)_{\g}$, this means that
the $\Sigma$-Hamiltonian flow of the  functions $\Phi_\I^z$
gives the
$G$-action and  the map
\begin{equation}\label{eq:554} 
\g\to R_{\I-hol}(\O),\qquad z\mapsto\Phi_\I^z
\end{equation}
is a complex Lie algebra homomorphism with respect to the
Poisson bracket on $R_{\I-hol}(\O)$ defined by $\Sigma$.
 
The space $\Mk$ has a natural $U$-invariant hyperkaehler
structure $(g,\J_1,\J_2,\J_3,\om_1,\om_2,\om_3)$; see
\cite[Remark 2, pg 476]{Kr2}  and \cite{H}.  Kronheimer
discovered
\begin{theorem}\label{thm:56} 
\cite{Kr2}  
\item{\rm(i)} The map $\ze$ in {\rm(\ref{eq:533})} is an 
$(SO(3)\times U)$-equivariant  smooth embedding of  
manifolds.
\item{\rm(ii)}
The three maps $\ze_1,\ze_2,\ze_3$ are the moment maps for 
the $U$-action with respect to the three Kaehler forms
$\om_1,\om_2,\om_3$ on $\Mk$.
\item{\rm(iii)} 
For $a=1,2,3$,  the holomorphic moment 
map $\Phi_a:\Mk\to\g$ given by  {\rm(\ref{eq:256})} is $1$-to-$1$ and 
has image equal to $\O$. Thus we get a   $G$-equivariant
holomorphic symplectic  isomorphism
\begin{equation}\label{eq:561} 
\Phi_a:\Mk\to\O
\end{equation}
from $(\Mk,\J_a,\Om_a)$    to $(\O,\I,\Sigma)$.
Here $G$ acts on $\Mk$ by the $\J_a$-complexification of the
$U$-action.
\item{\rm(iv)} 
The $SO(3)$-action 
\begin{equation}\label{eq:562} 
SO(3)\to\Diff\Mk 
\end{equation}
preserves the Riemannian  metric $g$ and induces the standard   
transitive action of $SO(3)$ on the  $2$-sphere $\cS_{\Mk}$ of
Kaehler  structures on $\Mk$.
\end{theorem}

In fact, (ii) determines uniquely the $U$-invariant 
hyperkaehler structure on $\Mk$. Let $\cC_a=\cC_{\bj_a}$.     
We further note

\begin{corollary}\label{cor:56} 
\item{\rm(i)} The  map $\ze$ is $\RP$-equivariant with respect
to {\rm(\ref{eq:541})} and the Euler scaling action on
$\bi\u\oplus\bj\u\oplus\bk\u$, i.e.,
\begin{equation}\label{eq:563} 
\ze(\la\diamond A)=\la\ze(A)
\end{equation}
\item{\rm(ii)} The
isomorphism $\Phi_a$  intertwines the   product action of 
$\RP\times\cC_a$ on $\Mk$ with the Euler scaling  action of 
$\C^*$  on $\O$.
\item{\rm(iii)} The $SO(3)$-action on $\Mk$ is free. 
\end{corollary}
\begin{proof} 
(i) and (ii) are routine to verify. The Euler $\C^*$-action on
$\O$ has only trivial isotropy groups, and so (ii) implies that
the action of $\cC_a$, and hence of all of $SO(3)$, has only
trivial isotropy groups. A  compact group action with only
trivial isotropy groups   is necessarily free.
\end{proof}
 
 Next we define a left action of $\H^*$ on $\Mk$ by
\begin{equation}\label{eq:571} 
h\cdot A_w=|h|^2\diamond A_{h\i{w}h} 
\end{equation}
This is the product of  the  square of the $\RP$-action
(\ref{eq:541}) with the $SU(2)$-action defined by the $SO(3)$-action 
on $\Mk$  the spin homomorphism
\begin{equation}\label{eq:572} 
SU(2)\to SO(3) 
\end{equation}
We also have an $\H^*$-action on $\bi\u\oplus\bj\u\oplus\bk\u$,
$(h,p)\mapsto h\bullet p$ defined by (\ref{eq:451}).
 
Let $\eta$ be the infinitesimal generator of the  square of the
action (\ref{eq:541}).

\begin{corollary}\label{cor:57} 
 The left $\H^*$-action on $\Mk$ defined in {\rm(\ref{eq:571})} satisfies
\begin{equation}\label{eq:573} 
\ze(h\cdot A)=h\bullet\ze(A) 
\end{equation}
commutes with the action of $U$, and gives $\Mk$   the
structure of a hyperkaehler cone. 
\end{corollary}
\begin{proof}
The $\H^*$-action on $\Mk$  clearly commutes with the 
$U$-action and satisfies (\ref{eq:532}) because of the
$SO(3)$-equivariance in Theorem \ref{thm:56}(i) and the  
$\RP$-equivariance in Corollary \ref{cor:56}(i).

Next we check the hyperkaehler cone axioms (\ref{eq:342}).
By Theorem \ref{thm:56}, the  $SU(2)$-action transforms the 
tensors $\J\qq$, $\om\qq$, $g$ according to (\ref{eq:342}).  
The action of the subgroup $\C^*_a=\C^*_{\bj_a}$
of $\H^*$  identifies (because of  Corollary \ref{cor:56}(ii)) 
with the 
\emph{square} of  the  Euler $\C^*$-action on $\O$ under  
$\Phi_a$. Hence the action of $\C^*_a$
is $\J_a$-holomorphic. Furthermore, this implies that
the action of the $\RP$ subgroup   of $\H^*$  
preserves $\J_1,\J_2,\J_3$,  and  transforms 
$\Om_1,\Om_2,\Om_3$, and hence 
$\om_1,\om_2,\om_3$, by the degree $2$ character.  
Thus the $\RP$-action also satisfies the axioms (\ref{eq:342}).
\end{proof}

Let $\eta$ be the infinitesimal generator of the action on 
$\Mk$ of the $\RP$-subgroup of $\H^*$, i.e., of the square of 
the action (\ref{eq:541}).   We will say a function  $f$ on 
$\Mk$ is  \emph{homogeneous of degree $r$} if this 
$\RP$-action transforms 
$f$ by the  degree $r$ character, i.e.,  $\eta f=rf$.

\begin{corollary}\label{cor:58} 
 The hyperkaehler manifold $\Mk$ admits a 
$U$-invariant hyperkaehler potential $\rho$, unique up to
addition of a constant. The further condition that $\rho$ is
homogeneous of degree $2$ determines $\rho$ uniquely. 
\end{corollary}
\begin{proof}
This follows by   Corollary \ref{cor:57} and Proposition 
\ref{prop:46} as soon as we check that
\[H^{1}_{DeRham} \Mk =0\] 
But $\Mk$ is diffeomorphic to $\O$ by Theorem 
\ref{thm:56}(iii) and we have the  well-known fact 
\begin{lemma}\label{lem:58} 
The fundamental group of a complex
nilpotent orbit $\O$ is finite. Consequently 
$H^{1}_{DeRham} \O =0$. 
\end{lemma}
This finishes the proof of Corollary \ref{cor:58}.
\end{proof}

\section{Complex Conjugation on the Instanton Space
$\Mk$}\label{sec:6}
\setcounter{equation}{0}

We will introduce a family of complex conjugation
maps on $\Mk$.

We start with the non-trivial Lie algebra involution
$\vsig$ of $\so(3)$ given by
\begin{equation}\label{eq:621} 
\vsig(a\bj_1+b\bj_2+c\bj_3)=-a\bj_1+b\bj_2-c\bj_3
\end{equation} 
(The exact choice is immaterial, as any two choices are
$SO(3)$-conjugate.) We assume now that the induced 
involution on $\ka(\so(3))$ extends to a Lie algebra involution 
$\vth$ of $\u$. Then $\vth$ determines a splitting
\begin{equation}\label{eq:622} 
\u=\a\oplus\b 
\end{equation}
where $\a$ is the $1$-eigenspace of $\vth$ and $\b$ is the
$(-1)$-eigenspace of $\vth$. This gives the bracket relations 
\begin{equation}\label{eq:623} 
[\a,\a]\subset\a,\qquad [\a,\b]\subset\b,\qquad 
[\b,\b]\subset\a 
\end{equation}

\begin{lemma}\label{lem:62} 
The  involution $\vth$ of $\u$ induces a
hyperkaehler involution $\Th$ of $\Mk$ which commutes with
the $SO(3)$-action.
\end{lemma}
\begin{proof}  It follows from the instanton differential equations 
(\ref{eq:523}) and the boundary conditions 
(\ref{eq:525}) that   any Lie algebra
automorphism of $\u$ induces a hyperkaehler diffeomorphism 
of $\Mk$. Clearly $\Th$ commutes with the $SO(3)$-action.
\end{proof}

Next we  extend  $\vth$ in a $\C$-antilinear manner to get a 
real Lie  algebra involution $\nu$ of $\g$. We have
\begin{equation}\label{eq:624} 
\g=\u\oplus i\u=\a\oplus\b\oplus i\a\oplus i\b 
\end{equation}
and  so for $x,x'\in\a$ and $y,y'\in\b$ we get
\begin{equation}\label{eq:625} 
\nu(x+y+ix'+iy')=x-y-ix'+iy' 
\end{equation}
 
The point is that $\nu$ is a   complex conjugation map
on the Lie algebra $\g$. The corresponding real form of 
$\g$ consisting of $\nu$-fixed vectors is
\begin{equation}\label{eq:626} 
\gR=\a\oplus i\b 
\end{equation}
Now  $\gR$ is a non-compact real form of $\g$, unlike   $\u$ 
which is a compact  real form. Indeed    the Killing form of $\g$,
being negative-definite on $\u$, is  clearly indefinite on $\gR$.
From now on, we call $\nu$ ``complex conjugation" and we
put  $\oz=\nu(z)$ for $z\in\g$.

Clearly $\gR$ is compatible with $\ka$ in the sense that
$\ka(\so(3))$ is complex conjugation stable in $\u$; in fact
$\ka(\vsig(w))=\vth(\ka(w))$ for $w\in\so(3)$.  Hence complex
conjugation preserves $\O$ inside $\g$. We conclude
\begin{lemma}\label{lem:63} 
Our map $\nu:\g\to\g$ induces an
antiholomorphic diffeomorphism  of $\O$, which we   again 
call $\nu$. So $\nu:\O\to\O$ is a complex conjugation map.
\end{lemma}

See  \S\ref{sec:7} and (\ref{eq:731})  for how this arises in
practice. We work with (\ref{eq:626})
 now because the essential symmetry of
the picture  is more apparent this way (just as for the 
Sekiguchi correspondence).

For $a=1,2,3$, we call pull back complex conjugation
on $\O$ through the holomorphic isomorphism $\Phi_a$
in (\ref{eq:561}) to get an involution
\begin{equation}\label{eq:641} 
\ga_a:\Mk\to\Mk,\qquad  \ga_a=\Phi_a^*\nu 
\end{equation}
So $\Phi_a(\ga_aA)=\nu(\Phi_aA)$.
Then $\ga_a$ is $\J_a$-antiholomorphic and so
$\ga_a$ defines a complex conjugation map on $\Mk$ with
respect to $\J_a$.

Our aim is to figure out how the hyperkaehler data transforms
under the involutions $\ga_a$.  To this end, we find a formula
for $\ga_a$ in terms of $\Th$ and the $SO(3)$-action on $\Mk$.
In particular,  the rotation in $\Him$ about the $\bj_a$-axis
defines via (\ref{eq:562}) an  isometric involution
\[R_a:\Mk\to\Mk\]

\begin{proposition}\label{prop:64} 
Let $(abc)$ be a cyclic permutation of $1,2,3$. Then
\begin{equation}\label{eq:642} 
\ga_a=\Th R_b
\end{equation}
\end{proposition}
\begin{proof}
Let $\vth:\g\to\g$ be the  involution defined by $\C$-linearly
extending $\vth:\u\to\u$. Then  $\vth$ induces a biholomorphic
diffeomorphism of $\O$, which we again call $\vth$. It is easy 
to see that 
\begin{equation}\label{eq:643} 
\Th=\Phi_a^*\vth
\end{equation}
The composition 
$\be=\nu\vth$ is the $\C$-antilinear involution  $\be:\g\to\g$ 
with fixed space $\u$, i.e, $\be$ is complex conjugation with
respect to $\u$. So $\be$ induces an antiholomorphic
involution $\be$ of $\O$. We now see that (\ref{eq:642}) holds if and
only if  
\begin{equation}\label{eq:644} 
R_b=\Phi_a^*\be 
\end{equation}

To compute $\Phi_a^*\be$, we first write out the map 
$\Phi_a$ according to (\ref{eq:255}):
\begin{equation}\label{eq:645} 
\Phi_a=\ze_b+i\ze_c
\end{equation}
It follows that  the  diffeomorphism $\tbe=\Phi_a^*\be$
of  the instanton space $\Mk$ satisfies
\[\tbe^*\ze_b=\ze_b\AND \tbe^*\ze_c=-\ze_c\]But the  rotation $R_b$
transforms $\ze_b$ and $\ze_c$ in exactly the same way.  This proves
(\ref{eq:644}).
\end{proof}
  
Our  formula (\ref{eq:642}) immediately implies, because of 
Theorem \ref{thm:56}(iv):
 
\begin{corollary}\label{cor:65} 
Let $(abc)$ be a cyclic permutation of $1,2,3$.
Then  the involution $\ga_a$ is an isometry of $\Mk$ which
preserves $(\J_b,\om_b)$ but  negates $(\J_a,\om_a)$ and 
$(\J_c,\om_c)$. Furthermore, $\ga_a$ commutes with the 
action of $\cC_b$.  

In particular   $\ga_a$ 
is a  complex conjugation map with respect to $\J_a$ which is 
also $g$-isometric, $\J_b$-holomorphic, $\om_b$-symplectic
and $\cC_b$-invariant.
\end{corollary}

Let $\Mk^{\ga_a}$ be the subset of $\Mk$ fixed by $\ga_a$.
 
\begin{corollary}\label{cor:66} 
\item{\rm(i)} $\Mk^{\ga_1}$ is a smooth submanifold of $\Mk$ 
with  finitely many connected components, and each
component has half the dimension of $\Mk$. 
\item{\rm(ii)} $\Mk^{\ga_1}$ is a symplectic real form of
$\Mk$ with  respect to $(\J_1,\Om_1)$.
\item{\rm(iii)}
$\Mk^{\ga_1}$ is a   Kaehler submanifold of $\Mk$ with    
respect to $(g,\J_2,\om_2)$.
\item{\rm(iv)} $\Mk^{\ga_1}$ is a complex Lagrangian
submanifold of $\Mk$ with  respect to
$(\J_2,\Om_2)$. 
\end{corollary}
\begin{proof}
If $\al:M\to M$ is a  finite-order automorphism of a smooth
manifold, then the fixed space $M^\al$ is a smooth submanifold.
The tangent space of $M^\al$ at $m$ is the subspace of
$\alpha$-fixed vectors  in the tangent space of $M$ at $m$.  

Now Corollary \ref{cor:65}  implies that $\Mk^{\ga_1}$ is a
Kaehler submanifold with respect to $(\J_2,\om_2)$. So
$\Mk^{\ga_1}$ is a $\J_2$-complex submanifold and  the
restriction of  $\om_2$ to it is a real symplectic form $\pi$.
Furthermore $\Mk^{\ga_1}$ is  a real form 
with respect to $\J_1$ and $\J_3$ and
$\om_1$ and $\om_3$ vanish on $\Mk^{\ga_1}$. So
$\Mk^{\ga_1}$ is Lagrangian with respect to 
$\Om_2=\om_3+i\om_1$ and  $\Om_1$ restricts to 
$\pi$. This proves (ii)-(iv).

Finally we complete the proof of (i). First,   $\Mk^{\ga_1}$,
being a real form of $\Mk$, has half the dimension. Next,
$\Mk^{\ga_1}$ has only  finitely many components since 
 $\ga_1=\Phi_1^*\nu$ by  (6.4.1),  and  so 
\begin{equation}\label{eq:661} 
\Phi_1(\Mk^{\ga_1})=\O\cap\gR 
\end{equation}
But  $\O\cap\gR$ is a real algebraic variety and so has
only  finitely many components.
\end{proof}

 Corollary \ref{cor:66}  above and the next result below encode
Vergne's result  in  ``hyperkaehler language".

\begin{proposition}\label{prop:67} 
We have
\begin{equation}\label{eq:671} 
\begin{array}{cl}
\Phi_1(\Mk^{\ga_1})&=\; \O\cap(\a\oplus i\b)\\[3pt] 
\Phi_2(\Mk^{\ga_1})&=\; \O\cap(\b\oplus i\b)\\[3pt]
\Phi_3(\Mk^{\ga_1})&=\; \O\cap(\b\oplus i\a)
\end{array}
\end{equation}
\end{proposition}

We already observed   in (\ref{eq:661}) that the first equality 
is clear, but the others reveal a subtle structure.

\begin{proof}
We start by computing how the functions $\ze_a^u$, $u\in\u$,
transform under $\ga=\ga_1$. For any function $f$ on $\Mk$, 
we   put $f^\ga=\ga^*f=f\circ\ga$.

By (\ref{eq:641}), we have $\ga=\ga_1=\Th R_{b}$. So
$\ga^*=R_{b}^*\Th^*$.
We find  $\Th^*\ze_a^u=\ze_a^{\vth u}$ and
$R_{b}^*\ze_a^u=t_a\ze_a^u$ where  $t_1=t_3=-1$ and   
$t_2=1$. Therefore, writing $u=x+y$ where $x\in\x$ and
$y\in\y$, we get
\begin{equation}\label{eq:672} 
\left(\ze^{x+y}_1\right)^\ga=\ze^{-x+y}_1,\qquad
\left(\ze^{x+y}_2\right)^\ga=\ze^{x-y}_2,\qquad
\left(\ze^{x+y}_3\right)^\ga=\ze^{-x+y}_3
\end{equation}
Let $u'=-x+y$ and $u''=x-y$. Then (\ref{eq:672}) says
\begin{equation}\label{eq:673} 
\ze_1(\ga A)=\ze_1(A)',\qquad
\ze_2(\ga A)=\ze_2(A)'',\qquad
\ze_3(\ga A)=\ze_3(A)',\qquad 
\end{equation}

Now we can compute the involutions $\al_b$, $b=1,2,3$, of 
$\O$ defined by $\Phi_a^*\al_a=\ga_1$. Of course $\al_1=\nu$
by the definition (\ref{eq:641}) of $ \ga_1$. Using   
(\ref{eq:255}) and (\ref{eq:673}) we  find
\begin{equation}\label{eq:674} 
\begin{array}{c}
\Phi_2(A)=\ze_3(A)+i\ze_1(A) \Rightarrow
\Phi_2(\ga A)=\ze_3(A)'+i\ze_1(A)'\\[6pt]
\Phi_3(A)=\ze_1(A)+i\ze_2(A) \Rightarrow
\Phi_3(\ga A)=\ze_1(A)'+i\ze_2(A)'' 
\end{array}
\end{equation}

Consequently,  the involutions $\al_2,\al_3$ of $\O$ are 
induced by the involutions $\al_2,\al_3$ of $\g$  given by,
for $x_1,x_2\in\a$ and $y_1,y_2\in\b$,   
\begin{equation}\label{eq:674b}  
\begin{array}{c}
\al_2(x_1+y_1+ix_2+iy_2)=-x_1+y_1-ix_2+iy_2\\[3pt]
\al_3(x_1+y_1+ix_2+iy_2)=-x_1+y_1+ix_2-iy_2
\end{array}
\end{equation}
So the fixed spaces are
\begin{equation}\label{eq:675} 
\g^{\al_1}=\a\oplus i\b, \qquad
\g^{\al_2}=\b\oplus i\b, \qquad
\g^{\al_3}=\b\oplus i\a 
\end{equation}
For $b=1,2,3$, we have 
\begin{equation}\label{eq:676} 
\Phi_b(\Mk^{\ga})=\O^{\al_b}=\O\cap\g^{\al_b}
\end{equation}
This completes the proof.
\end{proof}

We set
\begin{equation}\label{eq:681} 
P^1=\O\cap(\a\oplus i\b),\qquad P^2=\O\cap(\b\oplus i\b),
\qquad P^3=\O\cap(\b\oplus i\a) 
\end{equation}

Now combining Proposition \ref{prop:67} 
with  Corollary \ref{cor:66}, we get
\begin{corollary}\label{cor:68} 
For $a=1,2,3$, 
$P^a$ is the disjoint union of finitely many connected 
components. Each component  of  $P^a$ is a smooth real
submanifold of $\O$  of half the dimension.  Each component is
stable under the dilation action of $\RP$ on $\g$.

$P^1$ is a symplectic  real form  of $(\O,\I,\Sig)$ and $P^2$   
is a holomorphic  Lagrangian submanifold of $(\O,\I,\Sig)$.
\end{corollary}

Let $\Uvth$ be the subgroup of $U$ commuting with $\vth$;
then $\Uvth$ is the subgroup with Lie algebra equal to $\a$.

Let $\cs\in SO(3)$
be  the element which gives the cyclic  permutation
\begin{equation}\label{eq:691} 
\cs(a\bi+b\bj+c\bk)=c\bi+a\bj+b\bk 
\end{equation}
Then $\cs$ defines an automorphism $\cs:\Mk\to\Mk$ by 
(\ref{eq:562}).

\begin{corollary}\label{cor:69} 
\cite{Ve}
We have   $(\RP\times \Uvth)$-equivariant diffeomorphisms
\begin{equation}\label{eq:692} 
\Phi_b\Phi_a\i:P^a\to P^b 
\end{equation}
This sets up a bijection between the
connected components of $P^a$ and those of $P^b$.

We get the following commutative diagram of
diffeomorphisms where $(abc)$ is a cyclic permutation of
$1,2,3$:
\begin{equation}\label{eq:693} 
\begin{array}{ccc}
\Mk^{\ga_a}&\mapright{\cs\i}&\Mk^{\ga_b}\\[5pt]
\mapdown{\Phi_a}&&\mapdown{\Phi_a}\\[5pt]
P^a&\mapright{\Phi_b\Phi_a\i}&P^b 
\end{array}
\end{equation}
\end{corollary}
\begin{proof}We have
\begin{equation}\label{eq:694} 
\cs\cdot(A_1(t),A_2(t),A_3(t))
=(A_3(t),A_1(t),A_2(t)) 
\end{equation}
and so we get the commutative diagram
\begin{equation}\label{eq:695} 
\begin{array}{ccc}
\Mk&\mapright{\cs\i}&\Mk\\[5pt]
\mapdown{\Phi_a}&&\mapdown{\Phi_a}\\[5pt]
\O&\mapright{\Phi_b\Phi_a\i}&\O
\end{array}
\end{equation}
Now (\ref{eq:693}) follows by restriction.
\end{proof}

\section{Real Nilpotent Orbits and the Vergne
Diffeomorphism}\label{sec:7}
\setcounter{equation}{0}
 
 From now on,  we  single out  in Theorem \ref{thm:56}(iii)  the
holomorphic symplectic isomorphism
\begin{equation}\label{eq:721} 
\Phi_1:\Mk\to\O 
\end{equation}
Now by means of  $\Phi_1$  we transfer the hyperkaehler
structure on $\Mk$  over to $\O$. Thus we get   a  set
\begin{equation}\label{eq:722} 
(g,\I,\J,\K,\omI,\omJ,\omK) 
\end{equation}
of hyperkaehler data on $\O$ with a corresponding  
$2$-sphere $\cS$ of complex structures on $\O$.
We call this the \emph{instanton  hyperkaehler  structure}
on $\O$.  

By Theorem \ref{thm:56},  $\I=(\Phi_1)_*\J_1$ is the natural complex
structure  on  $\O$ discussed in \S\ref{sec:5} and the holomorphic
symplectic KKS form is 
\begin{equation}\label{eq:723} 
\Sigma=\omJ+i\omK 
\end{equation}
The three moment maps $\ze_1,\ze_2,\ze_3:\Mk\to\u$ 
transfer by means of  $\Phi_1$ to three $U$-equivariant
moment maps
\begin{equation}\label{eq:724} 
\ze_{\I}:\O\to\u,\qquad \ze_{\J}:\O\to\u,\qquad\ze_{\K}:\O\to\u   
\end{equation}
We recover the   natural embedding $\Psi:\O\to\g$ from two
of these maps as
\begin{equation}\label{eq:725} 
\Psi=\ze_{\J}+i\ze_{\K}   
\end{equation}
However the third moment map $\ze_{\I}$ is a mystery.

The $SO(3)$-action (\ref{eq:562}) on $\Mk$ transfers by means of 
$\Phi_1$ to a smooth action on $\O$ 
\begin{equation}\label{eq:726} 
SO(3)\to\Diff \O
\end{equation}
which commutes with the $U$-action.
We call this the\emph{ Kronheimer $SO(3)$-action} on $\O$. 

We emphasize that this $SO(3)$-action depends on a choice
of maximal compact subgroup of $U$ of $G$.  The Kronheimer 
$SO(3)$-action remains quite mysterious.

We   can rewrite (\ref{eq:626})  as
\begin{equation}\label{eq:731} 
\gR=\kR\oplus\pR 
\end{equation}
where $\kR=\a$ and $\pR=i\b$. 
This is a Cartan decomposition of the real semisimple Lie 
algebra $\gR$ in the usual sense:
$\kR$ is  a maximal compact Lie subalgebra of $\gR$ and
$\pR$ is an $\ad\kR$-stable complementary subspace.
Then (\ref{eq:622}) becomes
\begin{equation}\label{eq:732} 
\u=\kR\oplus i\pR 
\end{equation}
In practice, one may  start  from (7.3.1) and then constructs 
$\u$ by (\ref{eq:732}).    Then $\kR$ and $\pR$
are real forms of $\k=\a\oplus i\a$ and $\p=\b\oplus i\b$ 
and the complexification of (\ref{eq:731})  is
\begin{equation}\label{eq:733} 
\g=\k\oplus\p
\end{equation}
   
Let $(\,,\,)_{\gR}$ be the Killing form of $\gR$; this coincides
with the restriction to $\gR$ of the complex Killing form 
$(\,,\,)_{\g}$.

We know that $\Phi_1$ identifies $\Mk^{\ga_1}$ with
$\O\cap\gR$ by (\ref{eq:661}) or Proposition  \ref{prop:67}. 
Let  $\OR$ be a  connected component of  $\O$. Then $\OR$ is
stable under the  dilation action of $\RP$ on $\gR$. 
Moreover $\OR$ is an orbit under the adjoint action of $\GR$.
This follows easily since 
$\OR$ is a connected real form of $\O$  and $\O$ is an adjoint
orbit of $G$.   

Hence $\OR$ is a real nilpotent orbit of  $\GR$.  
Every real nilpotent orbit arises in this way. Let $\sig$ be the
real $\GR$-invariant KKS symplectic form  on $\OR$. Let 
\begin{equation}\label{eq:741} 
\phi:\OR\to\gR
\end{equation}
be the natural embedding. Then $\sig$ is the unique
symplectic    form   on $\OR$ such that the adjoint
action of $\GR$ on $\OR$ is  Hamiltonian with
moment map $\phi$. In terms of   the  component  function 
$\phi^w$, $w\in\gR$, defined by $\phi^w(z)=(w,z)_{\gR}$, this
means that the $\sig$-Hamiltonian flow of the
functions $\phi^w$ gives the $\GR$-action and  the map
\begin{equation}\label{eq:742} 
\gR\to\Cinf(\OR),\qquad w\mapsto\phi^w
\end{equation}
is a  Lie algebra homomorphism with respect to the
Poisson bracket on $\Cinf(\OR)$ defined by $\sig$.
We easily see that
\begin{lemma}\label{lem:74} 
The real nilpotent orbit $(\OR,\sig)$ is a symplectic
real form of the complex nilpotent orbit $(\O,\I,\Sig)$. I.e.,
$\OR$ is a real form of $\O$ and  the real part of $\Sig$  restricts
to $\sig$ on $\OR$ while the imaginary part of $\Sig$ vanishes
on $\OR$. Thus
\begin{equation}\label{eq:743} 
\sig=\Sig|_{\OR}
\end{equation}
\end{lemma}

\begin{theorem}\label{thm:75} 
$\OR$ is a complex submanifold of
$\O$ with respect to the complex structure $\J$ on $\O$.
Consequently the Kaehler structure $(g,\J,\omJ)$ on $\O$
defines by  restriction a $\KR$-invariant Kaehler structure on
$\OR$. 
 
The Kaehler form $\omJ|_{\OR}$ 
coincides with the  KKS symplectic form $\sig$.  
\end{theorem}
We call this Kaehler structure on $\OR$ the
\emph{instanton Kaehler structure}.

\begin{proof}
By (\ref{eq:743}) and (\ref{eq:723}), we have
\[\sig=\Re\Sig|_{\OR}=\Re(\om_\J+i\om_\K)|_{\OR}=
\om_\J|_{\OR}\]
The rest is an immediate consequence  of  our work in 
\S\ref{sec:6}.  Indeed, by Corollary \ref{cor:66}(iv),
$\Mk^{\ga_1}$ is a
$\J_2$-holomorphic complex submanifold of $\Mk$.
We have $\Phi_1^*\J=\J_2$ by the definition of $\J$. So
$\Phi_1(\Mk^{\ga_1})$ is a  $\J$-complex
submanifold of $\O$.
But we saw $\Phi_1(\Mk^{\ga_1})=\O\cap\gR$ in (\ref{eq:661}).
Thus each connected component of  $\O\cap\gR$ is a 
$\J$-complex submanifold of $\O$.
\end{proof}
\begin{remark}\label{rem:75}\rm 
Here is another proof that $\O\cap\gR$ is a 
$\J$-complex submanifold of $\O$. We can transport the
diffeomorphism $\cs:\Mk\to\Mk$  to a diffeomorphism
$\cs:\O\to\O$ by means of $\Phi_1$. Then Proposition 
\ref{prop:67} and   Corollary \ref{cor:69} imply
that $\O\cap\gR=\cs(\O\cap\p)$.  But $\O\cap\p$ is an
$\I$-complex submanifold of $\O$. So $\O\cap\gR$ is a
complex submanifold with respect to $\cs(\I)=\J$.
\end{remark}

\begin{example}\label{ex:75} \rm
We illustrate Theorem \ref{thm:75} in the context of  the
example  of flat hyperkaehler space 
discussed in Examples \ref{ex:23}, \ref{ex:26}, and \ref{ex:44}.
The case
$n=1$ is sufficient to show how this works.
 
Let  $\O\subset\sl(2,\C)$ be  the non-zero nilpotent
orbit of $SL(2,\C)$. We have a   covering map 
\[\pi:\R^4-\{0\}\to\O,\qquad \pi(a+b\bi+c\bj+d\bk)=
\left(\begin{array}{cc} uv&-u^2\\ v^2&-uv
\end{array}\right)\] where  $u=a+bi$ and $v=c+di$.
Let $\OR$ be the Euclidean connected component of
$\O\cap\sl(2,\R)$ containing $\pi(1+\bj)$.
For  $q=a+b\bi+c\bj+d\bk$ we find that
\[\pi(q)\in\OR\iff u^2,uv,v^2\in\R\mbox{ and } u^2+v^2>0  \iff 
u,v\in\R\iff b=d=0\] Thus
\[\pi\i(\OR)=(\R+\R j)-\{0\}\] 
Plainly,  $(\R+\R j)-\{0\}$ is a complex submanifold
of $\H-\{0\}$  with respect to  $\pm\J$ (and 
only $\pm\J$).  But $\pi$ is a covering of hyperkaehler
manifolds, and so it follows that $\OR$ is a complex submanifold
of $\H-\{0\}$  with respect to  $\pm\J$ (and 
only $\pm\J$).
\end{example}

We have a source of  $\J$-holomorphic functions on $\O$,
namely the component functions $\Phi_{\J}^z$,  $z\in\g$, 
of the $\J$-holomorphic moment map $\Phi_\J:\O\to\g$.
We next examine how these restrict to $\OR$.
For this, we  recall the complex Cartan decomposition
(\ref{eq:733}).
\begin{theorem}\label{thm:76} 
The functions $\Phi_{\J}^x$, $x\in\k$, vanish
identically  on  $\OR$. However for $v\in\p$, the  functions 
\begin{equation}\label{eq:761} 
f^v=i\Phi_{\J}^v|_{\OR} 
\end{equation}  separate the
points of $\OR$ and the corresponding holomorphic map
$\V:\OR\to\p$ defined by the $f^v$, so
that $(v,\V(w))_\g=f^v(w)$, is a locally closed embedding. 
 
The image $Y=\V(\OR)$ is a single
$K$-orbit in $\p$. $Y$ is a connected component of 
$\O\cap\p$ and so in particular $Y$ is stable under the dilation 
action of $\C^*$.   The resulting $(\KR\times\RP)$-equivariant
diffeomorphism
\begin{equation}\label{eq:762} 
\V:\OR\to Y
\end{equation}
is the    Vergne diffeomorphism.  
\end{theorem}
\begin{proof}
We know by Proposition \ref{prop:67} that $\Phi_\J$ gives a
diffeomorphism of $\O\cap\gR$  onto $\O\cap\p$. This means
that $\Phi_\J^x$, $x\in\k$,  vanishes on $\O\cap\gR$ while the
functions $\Phi_\J^v$, $v\in\p$,  embed   $\O\cap\gR$ into    
$\p$ as $\O\cap\p$.  But then the functions $i\Phi_\J^v$ 
do  the same thing. The diffeomorphism
\begin{equation}\label{eq:763} 
\V:\O\cap\gR\to\O\cap\p 
\end{equation}
defined by $i\Phi_\J$ is the Vergne diffeomorphism
discovered  in \cite{Ve}.

Now $Y=\V(\OR)$ is a connected component of $\O\cap\p$.
We know  by a well-known argument  that  $Y$ is
$K$-homogeneous. Indeed, $K$ acts on $\O\cap\p$ and at any
point $e\in Y$, we have a complex Lagrangian decomposition
\begin{equation}\label{eq:764} 
T_e\O=[\g,e]=[\k,e]\oplus[\p,e] 
\end{equation}
Also $[\k,e]\subset\p$ while $[\p,e]\subset\k$.  It follows
that $T_eY=[\k,e]$. So the $K$-action on $Y$ is  infinitesimally
transitive and hence transitive.
\end{proof}

\begin{remark}\label{rem:76}\rm  
 (i) We introduced the factor $i$ in (\ref{eq:761}) for
convenience (later in Theorems \ref{thm:78} and 
\ref{thm:94}) Essentially, the
factor $i$ arises because it is $i\OR$, not $\OR$, which most
naturally corresponds to $Y$ in the Sekiguchi 
correspondence.
 
(ii) The Vergne diffeomorphism in (\ref{eq:763}) recovers 
and explains the Kostant-Sekiguchi correspondence
\cite{sek}. 
\end{remark}

\begin{corollary}\label{cor:76} 
The $\KR$-action on $\OR$ complexifies with respect 
to $\J$ to give a transitive holomorphic action of $K$ on 
$\OR$. Then the Vergne diffeomorphism {\rm(\ref{eq:762})} is
$K$-equivariant. 
\end{corollary} 

 Theorem \ref{thm:76} and Proposition  \ref{prop:67} give   
\begin{corollary}\label{cor:77} \cite{Ve}
Let $A\in\Mk$. Then
\begin{equation}\label{eq:771} 
\Phi_1(A)\in\OR\quad\iff\quad\Phi_2(A)\in Y\quad\iff\quad
\Phi_3(A)\in i\OR
\end{equation}
\end{corollary}

Let
\begin{equation}\label{eq:781} 
\mu:\OR\to\kR
\end{equation}
be the projection map defined by the Cartan decomposition
(\ref{eq:731}). Then $\mu$ is a $\KR$-equivariant moment map for the
$\KR$-action (inside the $\GR$-action) on $\O$. 
\begin{theorem}\label{thm:78} Triple Sum Formula.
Let $w\in\OR$. Then the
Vergne diffeomorphism {\rm(\ref{eq:762})} satisfies
\begin{equation}\label{eq:782} 
w=\mu(w)+\half\V(w)+\half\ovl{\V(w)}
\end{equation}
\end{theorem}
In (\ref{eq:782}), we are taking the sum inside $\g$ of  the 
three vectors. In view of (\ref{eq:731}),  (\ref{eq:782}) says
that $\half\V(w)+\half\ovl{\V(w)}$ is the projection of $w$ to 
$\pR$.

\begin{proof}
Given $w\in\OR$, let $A\in\Mk$ be the (unique) instanton 
such that 
\begin{equation}\label{eq:783} 
w=\Phi_1(A)=\ze_2(A)+i\ze_3(A) 
\end{equation}
By Proposition \ref{prop:67}, we have (taking $\a=\kR$, 
$\b=i\pR$)
\begin{equation}\label{eq:784} 
\ze_1(A)\in i\pR,\qquad\ze_2(A)\in\kR,\qquad\ze_3(A)\in i\pR
\end{equation}
 
So the projection of $w$ to $\kR$ is
\begin{equation}\label{eq:784b} 
\mu(w)=\ze_2(A) 
\end{equation}
Now
\begin{equation}\label{eq:785} 
\V(w)=i\Phi_2(A)=i(\ze_3(A)+i\ze_1(A))=-\ze_1(A)+i\ze_3(A)
\end{equation}
Since $\ze_1(A)\in i\pR$ is pure imaginary and 
$i\ze_3(A)\in\pR$ is real, we get
\begin{equation}\label{eq:786} 
\ovl{\V(w)}=\ze_1(A)+i\ze_3(A) 
\end{equation}
So
\begin{equation}\label{eq:787} 
\half\left(\V(w)+\ovl{\V(w)}\right)=i\ze_3(A) 
\end{equation}
Now (\ref{eq:782}) is immediate.   
\end{proof}

The real symmetric pair $(\gR,\kR)$ (or the complex 
symmetric pair $(\g,\k)$) is called \emph{Hermitian} if   the
corresponding real symmetric space $\GR/\KR$ has a
$\KR$-invariant Hermitian structure. This amounts to the
condition that we can find  $x_0\in\Cent\kR$ such that
the eigenvalues of $\ad x_0$ on $\p$ are $\pm i$. 
Then we get a $\KR$-stable splitting
\[\p=\p^+\oplus\p^-\]
where $\p^{\pm}$ is the $\pm i$-eigenspace.
Then $\dimC\p^+=\dimC\p^-$ and moreover
$\p^+$ and $\p^-$ are  mutually 
contragredient  $\KR$-representations.

\begin{lemma}\label{lem:79} 
Suppose $(\g,\k)$ is Hermitian. Then
the following two $\R$-linear maps are inverses:
\[\pR\mapright{\iota}\p^+, \quad
u\mapsto \half\left(u-i[x_0,u]\right)\AND
\p^+\to\pR,\quad v\mapsto v+\ov\]
\end{lemma}

The Addition Formula  recovers in the 
Hermitian case the fact      

\begin{corollary}\label{cor:79} \cite[Prop. 6]{Ve}
Suppose $(\g,\k)$ is Hermitian and $\OR$ is such that
$Y=\V(\OR)$  lies in  $\p^+$.   Then
the Vergne diffeomorphism $\V:\OR\to Y$ is given by the
composition 
\[\OR\hookrightarrow\gR\to\pR\mapright{\iota}\p^+\]
where the middle map is the projection defined by
the Cartan decomposition {\rm(\ref{eq:731})}.
\end{corollary}
   
\section{ The KV Circle Action on $\OR$ and the Kaehler
Potential $\rhoo$}\label{sec:8}
\setcounter{equation}{0}

The  instanton Kaehler  structure $(\J,\sig)$ 
We let $\d=\del+\odel$ be the canonical splitting of $\d$ into
its $\J$-holomorphic and $\J$-antiholomorphic parts.
(So $\del=\del_\bj$ and $\odel=\odel_\bj$.) 

Corollary \ref{cor:58} and Corollary \ref{cor:56}(ii) give
\begin{corollary}\label{cor:82} 
The complex nilpotent orbit $\O$, equipped with its
$U$-invariant instanton  hyperkaehler  structure, admits a
$U$-invariant hyperkaehler potential $\rho:\O\to\R$, unique up
to addition of a constant. The further condition that $\rho$ is
homogeneous of degree $1$, with respect to the Euler
scaling action of $\RP$ on $\O$,  determines $\rho$ uniquely.
\end{corollary}

 Let
\begin{equation}\label{eq:831} 
\rhoo:\OR\to\R
\end{equation}
be the smooth function obtained by restricting the
hyperkaehler potential $\rho$  from Corollary \ref{cor:82}.
Plainly Theorem \ref{thm:56} and Corollary \ref{cor:82} give 
\begin{corollary}\label{cor:83} 
$\rhoo$ is a global Kaehler potential on $\OR$ so that 
\begin{equation}\label{eq:832} 
\sig=i\del\odel\rhoo 
\end{equation}
The function  $\rhoo$ is $\KR$-invariant and 
homogeneous of degree $1$ with respect to the Euler scaling
action of $\RP$ on $\OR$.
\end{corollary}

Next we consider how the Kronheimer $SO(3)$-action
on $\O$ from \S\ref{sec:7} transforms the submanifold $\OR$. 

\begin{proposition}\label{prop:84} 
The circle subgroup $\cC_{\bj}$ of $SO(3)$  
preserves $\OR$ and so defines a smooth group action
\begin{equation}\label{eq:841} 
\Pi:S^1\to\Diff\OR
\end{equation} 
We  call this the KV \(Kronheimer-Vergne\) $S^1$-action, and
we let $\th$ denote the infinitesimal generator.

The KV $S^1$-action  on $\OR$ is  free, commutes with
the $\KR$-action and is  Kaehler.
Moreover the KV $S^1$-action is Hamiltonian and $\rhoo$ is 
an equivariant moment map. I.e.,
\begin{equation}\label{eq:842} 
\th\rbot\sig+\d\rhoo=0
\end{equation}
or, equivalently, The Hamiltonian flow of $\rhoo$ is the KV
$S^1$-action.
\end{proposition}
We write the KV $S^1$-action as, for $w\in\OR$,  
\begin{equation}\label{eq:843} 
\Pi({\rm e}^{it})(w)={\rm e}^{it}\star w
\end{equation}
\begin{proof}
The first statement is clear by, e.g., Corollary \ref{cor:65}. The 
properties  of the KV $S^1$-action are then immediate from 
the  corresponding properties of the Kronheimer 
$SO(3)$-action on $\O$ or $\Mk$. 
\end{proof}

The $\J$-holomorphic functions $\Phi_\J^z$,
$z\in\g$, on $\O$ transform  under the action of $\cC_{\bj}$ by
the degree $1$ character. Consequently Proposition 
\ref{prop:84} gives
\begin{corollary}\label{cor:85} 
For $v\in\p$, the $\J$-holomorphic functions $f^v$  
on $\OR$ satisfy 
\begin{equation}\label{eq:851} 
\{\rhoo,f^v\}=if^v
\end{equation}
\end{corollary}

 We get insight into the KV $S^1$-action by trying to
find it inside the action of $\KR$.
For simplicity of  exposition, we assume in the next result that
$\g$ is  simple. If $\g$ is simple, then either 

(I) $(\g,\k)$ is Hermitian and $\Cent\kR=\R x_0$  where
$x_0\neq 0$ was chosen in \S7.9, or 

(II) $(\g,\k)$ is non-Hermitian, $\Cent\kR=0$ and $\p$ is  
irreducible as a $\KR$-representation.  

\begin{proposition}\label{prop:86} 
 Suppose $\g$ is a simple complex Lie algebra. Then 
\item{\rm(I)} If $(\g,\k)$ is Hermitian and
$Y\subset{\p^{\pm}}$ then the KV $S^1$-action on $\OR$ is
given by the center  of $\KR$ and
\begin{equation}\label{eq:861} 
\rhoo=\phi^{\mp x_0} 
\end{equation}
\item{\rm(II)} Otherwise, the  KV $S^1$-action on $\OR$ lies
outside the action of $\GR$. 
\end{proposition}
\begin{proof}
To prove that $\rhoo=\phi^{cx_0}$ where $c$ is constant, it is
necessary and sufficient to show that
\begin{equation}\label{eq:862} 
\{\rhoo,f^v\}=\{\phi^{cx_0},f^v\} 
\end{equation}  for all
$v\in\p$. This is because the functions $f^v$, $v\in\p$, 
separate the points of
$\OR$ and both $\rhoo$ and $\phi^{x_0}$ are
$\RP$-homogeneous of degree $1$.

Let $v\in\p$. Then $\{\phi^{x_0},f^v\}=f^{[x_0,v]}$ and so
in the Hermitian case we get
\begin{equation}\label{eq:863} 
\{\phi^{x_0},f^v\}=\left\{  
\begin{array}{ll}
\phantom{-}if^v&\mbox{ if  $v\in\p^+$}\\ 
-if^v&\mbox{ if $v\in\p^-$}
\end{array}\right.       
\end{equation}
We want to compare this with   (\ref{eq:851}). 
If $Y\subset\p^+$, then $f^v$ vanishes on $Y$ for  all $v\in\p^+$. 
For $v\in\p^-$, (\ref{eq:862}) holds with  $c=-1$.
So then $\rhoo=\phi^{-x_0}$. Similarly
$Y\subset\p^-$ gives $\rhoo=\phi^{x_0}$. In either case, it
follows that $\Cent\KR$ gives the Hamiltonian flow of $\rhoo$,
and this is just the $KV$ circle action.

Now, continuing the Hermitian case, suppose   $Y$ fails to lie in
$\p^+$ or $\p^-$. Then we cannot find
$c\in\R$  such that (\ref{eq:862}) holds for all $v\in\p$. Since $\rhoo$ 
is $\KR$-invariant  and (up to scaling) $\phi^{x_0}$ is the only
$\KR$-invariant function in $\phi^*(\gR)$, we see that
$\rhoo$  lies outside of $\phi^*(\gR)$. If $(\g,\k)$ is
non-Hermitian, then $\phi^*(\gR)$ has no non-zero
$\KR$-invariants, and so certainly $\rhoo$  lies outside of
$\phi^*(\gR)$.  

The condition that $\rhoo$  lies outside of
$\phi^*(\gR)$ means exactly that the Hamiltonian flow of
$\rhoo$ lies outside the action of $\GR$.
\end{proof}
 
Using $\rhoo$, we  next  prove in Theorem   \ref{thm:87}
below that the Vergne diffeomorphism produces an
embedding of $\OR$ into $T^*Y$ which realizes $T^*Y$ as a
symplectic complexification of $\OR$.

Let
\begin{equation}\label{eq:871} 
s:\OR\;\;\mapright{\V\times\mu}\;\; 
Y\times\kR\;\;\mapright{}\;\; Y\times\k
\end{equation}
be the smooth map built out of the Vergne diffeomorphism,
the moment map $\mu$ and the obvious inclusion
$\kR\hookrightarrow\k$.    
Then $s$ is an   embedding of    $\OR$ into
$Y\times\k$, since  already $\V$ is $1$-to-$1$.

On the other hand, we have a natural holomorphic 
embedding
\begin{equation}\label{eq:872} 
T^*Y\to Y\times\k 
\end{equation}
of the  cotangent bundle $T^*Y$ of  $Y$. 
Indeed, differentiating the $K$-action on $Y$  we get an
infinitesimal holomorphic vector field action
\begin{equation}\label{eq:eta} 
\k\to\Vect^{hol}\;\OR,\qquad x\mapsto\eta^x
\end{equation}   
The induced $K$-action on $T^*Y$ is holomorphic
Hamiltonian with   momentum functions given by the
holomorphic symbols of the vector fields $\eta^x$.
Let $M:T^*Y\to\k$ be the corresponding 
holomorphic moment map. Then (\ref{eq:872}) is the product
of the  canonical projection $T^*Y\to Y$ with $M$. It
follows since the $K$-action on $\OR$ is transitive
that the product map is $1$-to-$1$.

\begin{theorem}\label{thm:87} 
The image of $s$ in $Y\times\k$  lies inside $T^*Y$.
Furthermore    the resulting $\KR$-equivariant map
\begin{equation}\label{eq:873} 
s:\OR\to T^*Y 
\end{equation}
embeds $(\OR,\sig)$ as a totally real symplectic
submanifold of  $(T^*Y,\Om)$ so  that
\begin{equation}\label{eq:874} 
s^*(\Re\Om)=\sig\AND s^*(\Im\Om)=0 
\end{equation} 
where $\Om$ is the canonical holomorphic symplectic form 
on $T^*Y$.
\end{theorem}
\begin{proof}
We may regard   the vector fields  $\eta^x$ on $Y$ 
as vector fields   on $\OR$ using the Vergne diffeomorphism.
Showing that the image of $s$ lies in $T^*Y$ amounts to
showing that $\mu:\OR\to\kR$ is given by a real smooth
$1$-form $\be$ on $\OR$ in the sense that for all $x\in\kR$ 
we have
\begin{equation}\label{eq:875} 
\<\be,\,\eta^x\>=\mu^x\AND \<\be,\J\eta^x\>=0
\end{equation}
This follows easily because, since $K$ acts transitively on $\OR$
by Corollary \ref{cor:76},  the holomorphic tangent spaces of 
$(\OR,\J)$ are spanned by the vector fields $\eta^x$,
$\J\eta^x$, $x\in\kR$.

Furthermore,  we find  a $1$-form $\be$ giving $\mu$, then
\begin{equation}\label{eq:876} 
s^*(\Re\Th)=\be
\end{equation}
where $\Th$ is the holomorphic canonical $1$-form on $T^*Y$.
Then $\d\Th=\Om$ and so (8.7.5) holds if and only if
$\d\be=\sig$, i.e.,  if and only if $\be$  is a symplectic potential.

So the problem is to produce a symplectic potential $\be$ on
$\OR$ satisfying (\ref{eq:875}). We claim that
\begin{equation}\label{eq:877} 
\be=-\frac{i}{2}(\del-\odel)\rhoo
\end{equation}
works. First $\be$ is a symplectic potential since
\begin{equation}\label{eq:878} 
\d\be=-\frac{i}{2}(\del+\odel)(\del-\odel)\rhoo
=i\del\odel\rhoo=\sig
\end{equation}
It follows that the functions $\<\be,\eta^x\>$ are momentum
functions for the $\KR$-action on $\OR$. Consequently, for 
each $x\in\kR$, $\<\be,\eta^x\>$ is equal to $\mu^x$ up to
addition  of a constant. But both functions $\<\be,\eta^x\>$
and $\mu^x$ are homogeneous of degree $1$ with respect 
to the Euler  $\RP$-action, and so  $\<\be,\eta^x\>=\mu^x$.
Finally, we find
\begin{equation}\label{eq:879} 
\<\be,\J\eta^x\>=\<\J\be,\eta^x\>=
\half\<(\del+\odel)\rhoo,\eta^x\>=
\half\<\d\rhoo,\eta^x\>=\half\eta^x(\rhoo)=0 
\end{equation}
since $\J\del=i\del$, $\J\odel=-i\odel$ and $\rhoo$ is
$\KR$-invariant. This proves (\ref{eq:874}) and it follows easily that
$\OR$ is a totally real submanifold of $T^*Y$.
\end{proof}

\begin{corollary}\label{cor:87} 
The embedding $s$ realizes $(T^*Y,\Om)$ as a symplectic
complexification of $(\OR,\sig)$.
\end{corollary}

\section{Toward Geometric Quantization of Real Nilpotent 
Orbits}\label{sec:9}
 \setcounter{equation}{0} 

We can interpret Theorem \ref{thm:75} as providing a
polarization on $\OR$ in the sense of Geometric Quantization.
\begin{corollary}\label{cor:92} 
$\OR$, equipped with its KKS symplectic form $\sig$, admits two
transverse $\KR$-invariant complex polarizations,
namely the holomorphic and anti-holomorphic 
tangent spaces defined by the instanton Kaehler structure.
\end{corollary}

\begin{remark}\label{rem:92}\rm 
We have two rather different extensions of the 
Kaehler $\KR$-action on $\OR$ to a 
transitive action of a larger group. These extensions are given
by the $\GR$-action and the $K$-action. But
neither of these larger actions is Kaehler.
Indeed the $\GR$-action is symplectic but not holomorphic,
while the $K$-action  is holomorphic but not symplectic.
\end{remark}

The instanton Kaehler structure on $\OR$ provides the first 
step in our quantization program for $\OR$. In the spirit of 
 Geometric Quantization (GQ), the quantization problem on
$\OR$ is to convert $\R$-valued smooth functions $\phi$ on
$\OR$ into Hermitian operators $\Q(\phi)$ on a Hilbert space 
$\cH$ of
$\J$-holomorphic half-forms on $\OR$. This conversion must
satisfy the Dirac axiom that the  Poisson bracket of functions 
goes over into the commutator  of operators so that
\begin{equation}\label{eq:931} 
\Q\left(\{\phi,\psi\}\right)=i[\Q(\phi),\Q(\psi)] 
\end{equation}

By the well-known  No-Go Theorem, this conversion cannot be
carried out for all smooth functions (or even for all polynomial
functions on a real symplectic vector space) consistently so
that (\ref{eq:931}) is satisfied. However, we expect (for various
reasons) that the Hamiltonian functions $\phi^w$, $w\in\gR$,  
in (\ref{eq:742}) \emph{can be} quantized consistently, modulo some
``isolated anomalous" cases.

Quantization of the functions  $\phi^w$, $w\in\gR$, already
would solve the Orbit Method problem in representation 
theory of   attaching   to $\OR$  an   irreducible  unitary
representation (or finitely many such representations)
of the universal cover of $\GR$. Indeed, the operators
$i\Q(\phi^x)$  define a Lie algebra representation
of $\gR$ on the space of  $\KR$-finite vectors in $\cH$.

The larger GQ problem is to construct  the
operators $\Q(\phi^x)$ and the  Hermitian positive definite
inner product on (some subspace) of holomorphic
half-forms that gives rise to $\cH$.  

The  $\KR$-invariance of the instanton Kaehler
structure means that the Hamiltonian flow of the
functions $\phi^x$, $x\in\kR$ preserves $\J$.
Thus the Lie derivative  of the Hamiltonian vector field 
$\xi_{\phi^x}=\eta^x$  gives us a  natural choice for the
quantization, namely
\begin{equation}\label{eq:932} 
\Q(\phi^x)=-i\L_{\eta^x},\qquad x\in\kR 
\end{equation}

But the Hamiltonian flow of the ``remaining" functions
$\phi^v$, $v\in\pR$, does not preserve $\J$. Our strategy for
quantizing these remaining functions is to ``decompose"
them in terms of holomorphic and anti-holomorphic
functions on $\OR$.  The reason for this is that we expect
that a holomorphic function $f$ quantizes to multiplication by
$g$, while  an anti-holomorphic function $\ovl{g}$ quantizes
to the adjoint of multiplication by $f$.
(Another Dirac axiom is that $\Q(f)$ and $\Q(\ovl{f})$ are 
adjoint.)  Here $\Q(f)=f$ and $\Q(\ovl{g})=\Q(g)^*$ are densely
defined operators on $\cH$.  

A decomposition of $\phi^v$ with respect to holomorphic and
anti-holomorphic functions should be of the form
$\phi^v=\sum f_p\ovl{g}_p$ where
$f_p$ and $g_p$ are holomorphic functions. Of course, 
the sum could well be infinite and  unwieldy to quantize.

The key point is that  the decomposition of $\phi^v$  is
remarkably simple:

\begin{theorem}\label{thm:94} 
 Let $v\in\pR$. Then the Hamiltonian function
$\phi^v:\OR\to\R$ is the real part of a $\J$-holomorphic 
function on $\OR$.  In fact
\begin{equation}\label{eq:941} 
\phi^v=\Re f^v=\half\left(f^v+\ovl{f^v}\right) 
\end{equation}
where $f^v$ is the $\J$-holomorphic function defined in
Theorem {\rm \ref{thm:76}}. Moreover  
\begin{equation}\label{eq:942} 
f^v=\phi^v-i\{\rhoo,\phi^v\} 
\end{equation}
\end{theorem}
\begin{proof}  
The first part  is a corollary of the Triple Sum Formula 
Theorem \ref{thm:78}. Indeed, let $w\in\OR$. Then using    
(\ref{eq:782})   we find
\begin{eqnarray*}\label{eq:943} 
2\phi^v(w)&=&2(v,w)_\g=(v,\V(w)+\ovl{\V(w)})_\g=
(v,\V(w))_\g+\ovl{(v,\V(w))_\g}\\
&=&f^v(w)+\ovl{f^v(w)} 
\end{eqnarray*}
So $\phi^v=\half(f^v+\ovl{f^v})$.

Next, we have $\{-i\rhoo,f^v\}=f^v$ by (\ref{eq:851}) and so 
$\{-i\rhoo,\ovl{f^v}\}=-\ovl{f^v}$   by complex conjugation.
Then (\ref{eq:941}) gives
\begin{equation}\label{eq:944} 
\{-i\rhoo,\phi^v\}=\half(f^v-\ovl{f^v}) 
\end{equation}
Now  adding (\ref{eq:941}) and (\ref{eq:944})   we get 
(\ref{eq:942}).
\end{proof}

We  may choose a set $z_1,\dots,z_n$
of local holomorphic coordinates on $\OR$. Then
$\phi^v$ is a real analytic  function in the $2n$
coordinates $z_1,\dots,z_n,\oz_1,\dots,\oz_n$
and  (\ref{eq:941}) says that
\[\phi^v(z_1,\dots,z_n,\oz_1,\dots,\oz_n)
=f_v(z_1,\dots,z_n)+\ovl{f_v}(\oz_1,\dots,\oz_n)\]

\begin{corollary}\label{cor:94} 
The functions  $\phi^v$, $v\in\pR$, on $\OR$ are 
pluriharmonic.
\end{corollary}

\begin{remark}\label{rem:94}\rm 
(i) Alternatively, we can prove (\ref{eq:941}) right 
after Theorem \ref{thm:75}  in the following way. We find for $v\in\pR$
\[\phi^v=\Re\Phi_\I^v=\ze_\K^{iv}\]
because   (\ref{eq:255}) gives
$\Phi_\I^v=-i\Phi_\I^{iv}=-i\ze_\J^{iv}+\ze_\K^{iv}$. 
But also we have the $\J$-holomorphic function
\[i\Phi_\J^v=\Phi_\J^{iv}=\ze_\K^{iv}+i\ze_\I^{iv}\]
So $\phi^v=\Re f^v$ where $f^v=i\Phi_\J^v$.
 
\noindent (ii)  The condition $\phi^v=\Re F^v$ determines a
holomorphic function $F^v$ uniquely up to addition of a 
constant. The additional condition that $F^v$ is
homogeneous then forces $F^v=f^v$.
\end{remark}

We   summarize the way the KV $S^1$-action produces the
complex structure $\J$ on $\OR$ in the next Corollary.  
Let 
\begin{equation}\label{eq:971} 
\v=\phi^*(\pR)=\{\phi^v\,|\, v\in\pR\} 
\end{equation}
so that $\v$ is the space of Hamiltonian functions on $\OR$
corresponding to $\pR$. 
 
\begin{corollary}\label{cor:97} 
 Let $\vs\subset\Cinf(\OR)$ be the subspace spanned 
by all the translates of  $\v$ under the KV $S^1$-action.
If $\g$ is simple then in Cases {\rm(I)} and {\rm(II)} of
Proposition \ref{prop:86}  we find\item{\rm(I)}  $\vs=\v$ and
so
$\v\simeq\pR$ as $\KR$-representations.
\item{\rm(II)} We have the direct sum 
\begin{equation}\label{eq:972} 
\vs=\v\oplus\{\rhoo,\v\} 
\end{equation}
and so  $\v\simeq\pR\oplus\pR$ as $\KR$-representations. 
 
In either case, $\vs$ decomposes under the $KV$ action of 
$S^1\simeq SO(2)$ into a direct sum of copies of the 
$2$-dimensional rotation representation so that
$\v\simeq\R^2\oplus\cdots\oplus\R^2$.
Consequently the complexification  $\vC$
splits into the  direct sum 
\begin{equation}\label{eq:973} 
\vC=\vC^+\oplus\vC^- 
\end{equation}
where
$\vC^{\pm}=\{\phi\in\v\,|\, e^{i\th}\star\phi=
e^{\pm i\th}\phi\}$.
In Case {\rm I}, $\vC^{\pm}$ identifies with $\p^\pm$.
 
In either case, $\vC^{+}$ and $\vC^{-}$ are
complex-conjugate $\KR$-stable spaces of complex-valued 
Poisson commuting functions on $\OR$.   We have
\begin{equation}\label{eq:974} 
 \vC^+=\{f^v\,|\, v\in\p\}\AND\vC^-=\{\ovl{f^v}\,|\, v\in\p\}
\end{equation}
 
Finally  $\J$ is the unique complex
structure on $\OR$ such that  the functions $f^v$, $v\in\p$,
are $\J$-holomorphic.  
\end{corollary}

\begin{remark}\label{rem:97}\rm 
Corollary \ref{cor:97} says in particular that
\[\vs=\v+\{\rhoo,\v\}+\{\rhoo,\{\rhoo,\v\}\}+\cdots\]
Moreover, we could take this as the definition of $\vs$, and
then all the assertions in Corollary \ref{cor:97} read the
same.  
\end{remark}

\appendix

\section{Appendix: Real Algebraic Varieties and Nash 
Manifolds}\label{sec:A}
 
\subsection{Introduction}  In this appendix we
present some basic notions from the theory of
real algebraic varieties. Some references are
\cite{BoCR}, \cite{BeR}, \cite{BoE}.

On many points, we
follow the treatment in \cite{BoCR}. We add
material about real algebraic variety structures
associated to  complex algebraic varieties. 
Our discussion of Nash manifolds and Nash functions
is more general.
The last subsection \S\ref{ssec:last}
explains the application of this theory to
real groups and their orbits.

We present this material  here for lack
of an appropriate reference in the literature.

\subsection{Real Algebraic Sets}
A subset  $V\subset\R^n$ is a (real) \emph{algebraic set} if
$V$ is the set of common zeroes of some finite set 
of  real polynomial functions  
$P\in\R[x_1,\dots,x_n]$.
Then a function  $f:V\to\R$ is called \emph{regular} 
if there exist polynomial functions 
$P,Q\in\R[x_1,\dots,x_n]$ such that $Q$ has no 
zeroes on $V$ and 
$f(x)=P(x)/Q(x)$ for all $x\in V$.
The set of regular functions on $V$
forms an $\R$-algebra 
which we will call $A(V)$. 

Let $P(V)\subset A(V)$ be the image of the natural 
ring homomorphism $\R[x_1,\dots,x_n]\to A(V)$. 
Then 
\[P(V)\simeq \R[x_1,\dots,x_n]/I(V)\] where 
$I(V)\subset P(V)$ is the ideal of polynomial functions 
vanishing  on  $V$.
Notice that $A(V)$ is algebraic over $P(V)$;
indeed, if $f=P/Q$ then $Qf-P=0$.

A map $\phi:V\to V'$, where $V'\subset\R^m$ is an 
algebraic set,  is \emph{regular} if each of the 
component functions
$\phi_1,\dots,\phi_m$ 
is regular. An equivalent condition is that the pullback of  
a regular function on $V'$ is regular on $V$; then
$\phi$ induces an algebra homomorphism
$\phi^*:A(V')\to A(V)$. Conversely, any algebra
homomorphism
$p:A(V')\to A(V)$ defines uniquely a regular map
$V\to V'$.

An isomorphism of algebraic sets is a bijective
biregular map. Isomorphisms $V\to V'$ are in natural
bijection with algebra isomorphisms
$A(V')\to A(V)$.

We may also use the term \emph{real algebraic}
in speaking of  regular functions and maps. 

\subsection{ Real Affine Algebraic Varieties}
The \emph{Zariski topology} on an algebraic set 
$V\subset\R^n$ is defined just as in the complex 
case, so that the Zariski closed sets in $V$ are
precisely the algebraic sets in $\R^n$ which lie in $V$.
This topology is not Hausdorff but it is Noetherian
and hence  quasi-compact (every open cover has a finite
subcover)  as the polynomial ring
$R[x_1,\dots,x_n]$ is Noetherian.

Every Zariski closed
set  is closed in the usual Euclidean topology on
$\R^n$ defined by the Euclidean metric, as 
polynomials are continuous.  We will refer to open
sets, closed sets, etc as ``Zariski" or ``Euclidean" to
distinguish the two topologies.

A topological space $M$ is called irreducible
if $M$ cannot be written as  the union of two 
closed subsets different from $M$. We say an
algebraic set is \emph{irreducible} if it is irreducible
in the Zariski topology.

A \emph{regular function} on 
a Zariski open set $U\subset V$ is one of the form
$P(x)/Q(x)$ where
$P,Q\in\R[x_1,\dots,x_n]$ and $Q$ is nowhere
vanishing on $U$.
The set of regular functions on $U$ is closed under
composition and forms an
$\R$-algebra which we will denote $A_V(U)$.  

The assignment 
$U\mapsto A_V(U)$ defines a sheaf $\A_V$ of  
$\R$-algebras on $V$ with respect to its Zariski
topology. In particular, if $U_1,\dots,U_m$ is a
finite Zariski open cover of a Zariski open set
$U\subset V$ and $f$ is a real-valued function
on $U$ such that
$f|_{U_i}=P_i/Q_i$ then we can find
$P,Q\in\R[x_1,\dots,x_n]$ such that 
$Q$ is nowhere vanishing on $U$ and $f=P/Q$. 
Indeed, as $U_i\subset V$ is open, the complement
$V-U_i$ is the zero-locus of a finite set of
polynomials; let $F_i$ be the sum of their squares.
Then $U_i=V\cap(F_i\neq 0)$ and the polynomials
$P=\sum_{i}^{m} P_iQ_iF_{i}^2$ and
$Q=\sum_{i}^{m} Q_{i}^2F_{i}^2$ satisfy
our requirement.

Then $(V,\A_V)$ is a
ringed space in the usual sheaf theory sense.

Now we can define an (abstract) \emph{real affine 
algebraic variety}: this is  a pair $(X,\A_X)$ where
$X$ is an irreducible  topological space,
$\A_X$ is a sheaf of  $\R$-algebras of
$\R$-valued functions on $X$  and there exists an isomorphism
of ringed spaces from $(X,\A_X)$ to
$(V,\A_V)$ for some (irreducible) real algebraic set  $V$.

If $S\subset X$ is closed then $S$ identifies with 
an algebraic set of $\R^n$ inside $V$.
In this way, if $S$ is irreducible,  $S$ acquires a
canonical real algebraic affine variety 
structure; we call the corresponding structure sheaf  $\A_{X,S}$.

\subsection{Real Algebraic Varieties}
In complete analogy with the complex case, real
algebraic varieties are obtained by gluing
together affine ones.

A \emph{real algebraic variety} is a a pair $(X,\A_X)$
where $X$ is a Noetherian  irreducible topological space and
$\A_X$ is a sheaf of $\R$-valued functions on $X$
satisfying this condition: there exists a finite open cover 
$\{U_i\}_{i\in I}$ of $X$ such that for each $i$,
the ringed space $(U_i, \A_X|_{U_i})$ is a
real affine algebraic  variety. Then
$\A_X$ is called the \emph{structure sheaf} of
$X$. The sections of $\A_X$ are the \emph{regular
functions} on $U$. The topology of $X$ is then called 
the \emph{Zariski topology}. 

In speaking of real algebraic varieties, we may omit 
the modifiers ``real" or ``algebraic" when the 
context is clear. (However, often we will be dealing 
with complex algebraic varieties or real analytic
manifolds at the same time.)

A regular mapping between 
varieties $(X,\A_X)$ and $(Y,\A_Y)$ is 
a Zariski continuous mapping $\phi:X\to Y$ such that
if $U\subset Y$ is open and $f\in\A_Y(U)$ then
$\phi^*f=f\circ\phi\in\A_X(\phi\i U)$.
An isomorphism is a bijective biregular map.
We often speak of $X$ as the variety and leave
implicit its structure sheaf $\A_X$.

An \emph{affine open} set $U$ of $X$ is then a
Zariski open set $U$ such that
$(U,\A_X|_{U})$ is an affine variety.

Let $X$ be a a  real algebraic   variety Then we have the 
following examples of  real algebraic  subvarieties   of $X$.

(i) Suppose $W$ is     Zariski open in $X$. Then $W$ is again a 
variety where we define $\A_W$ by restriction of the 
structure sheaf of $X$.

(ii) Suppose $S$ is   an irreducible   Zariski closed set in $X$. 
Then $S$  is again a variety where for each open affine set 
$U\subset X$  we have $\A_S(S\cap U)=\A_{X,S}(S\cap U)$. 
If $X$ is affine then so is $S$.

(iii) Now (i) and (ii) imply that any Zariski locally
closed irreducible subset $W$ of    $X$ is again a
variety.  For such a subvariety we may write $\A_X|_{W}$
for the induced structure sheaf $\A_{W}$.
A regular map $\phi:X\to Y$ of varieties
is a \emph{locally closed embedding} if
$\phi(X)$ is a locally closed subvariety of $Y$
and $\phi$ defines an isomorphism
$(X,\A_X)\to (\phi(X),\A_{Y}|_{\phi(X)})$.

\subsection{Real Structures and Real Forms}
\label{ssec:real}

 A complex algebraic variety
$Z$ is \emph{defined over $\R$} if  $Z$ is 
equipped with an involution 
\[\kappa:Z\to Z\]  
called complex conjugation, which satisfies  
the following:
$Z$ admits a cover by complex algebraic  
$\kappa$-stable affine open subsets $U$
such that 

(i)  if $f\in R(U)$ then the function
$\of$ defined by $\of(u)=\ovl{f(\kappa(u))}$, 
$u\in U$, also lies in $R(U)$, and 

(ii) the map $R(U)\to R(U)$, $f\mapsto\of$,
is an $\R$-algebra involution  of $R(U)$.
In other words,  the real subspace 
$\{f\in R(U)\,|\, f=\of\}$ is both a real form 
and a real subalgebra of $R(U)$.

We  call this collection of real forms, or 
$\kappa$ itself, a \emph{real structure} on $Z$.
We write $\oz=\kappa(z)$ for $z\in Z$.
Geometric objects on $Z$ such as functions,
vector fields and differential forms are
\emph{defined over $\R$}, or \emph{real},
if they are stable under complex conjugation.

If $f:Z\to Z'$ is a complex algebraic morphism
of complex varieties defined over $\R$, then
$f$ is \emph{defined over $\R$}, or \emph{real},
if $f$ commutes with complex conjugation.

Let $Z^\kappa$  be the set of  real, i.e., $\kappa$-fixed
points in $Z$. We put $\Zr=Z^{\kappa}$.
Suppose  $\Zr$ is non-empty and Zariski irreducible.
Then   
$\Zr $ has a natural  structure of
real algebraic variety and we call $\Zr$ a \emph{real
form} of $Z$.
 
To see this, we first define a Zariski topology on 
$\Zr$ by the collection of sets 
\[\cS=\{\Ur=U\cap\Zr\,|\, U\in\T\}\]  where $\T$ is the
collection of  $\kappa$-stable Zariski open subsets  of
$Z$. For each affine $W\in\cS$ we
define $\A_{\Zr}(W)$ to be the space of
quotients $P/Q$ where $P$ and $Q$ are  
regular functions on some $U\in\T$ with
$\Ur=W$, $P$ and $Q$ real  (i.e., $\kappa$-fixed)  and $Q$
nowhere vanishing on $W$. This data defines a unique sheaf  
$\A_{\Zr}$ of  $\R$-valued functions on
$\Zr$ and then $(\Zr,\A_{\Zr})$ is a
real algebraic variety.

Notice that $\Zr$ is affine if $Z$ is affine.
In fact suppose $Z\subset\C^n$ is 
defined by the vanishing of
$P_1,\cdots,P_m\in\C[z_1,\dots,z_n]$ and also $Z$ 
is complex conjugation stable. Then 
$\Zr$ is the zero-locus in $\R^n$ of the 
$2m$ real polynomial functions defined by the real 
and imaginary parts 
$\hbox{Re}(P_1),\hbox{Im}(P_1),\dots,
\hbox{Re}(P_m),\hbox{Im}(P_m)\in\R[z_1,\dots,z_n]$. 

Clearly every real affine algebraic variety 
is of the form $\Zr$ for some complex affine
algebraic  variety $Z$ defined over $\R$.

The process $Z\mapsto\Zr$ is compatible with 
the usual operations on varieties. For instance, if
$\phi:Z\to Z'$ is a regular map of complex algebraic
varieties defined over $\R$, 
then the induced map $\phir:\Zr\to\Zr'$ is a
regular map of regular algebraic varieties.

Real structures often arise in the 
following way. Suppose $V$ is a complex vector
space and $\VR$ is a real form of $V$ with
corresponding complex conjugation map
$\kappa:V\to V$.  Then $\kappa$ defines a real 
structure on  every $\kappa$-stable (locally closed)
complex algebraic subvariety $X$ of $V$.  

Suppose  a complex  algebraic group $H$ acts on  
$Z$ and $H$ is  defined over $\R$.
We say the \emph{$H$-action on $Z$ is defined over 
$\R$} if the action morphism $H\times Z\to Z$ is defined
over $\R$. This happens if and only if   for every $h\in
H$, the transformations of $Z$ defined by $h$ and
$\ovl{h}$ are complex conjugate.

\subsection{The Complex Conjugate of a Complex
Variety}

Given a complex algebraic variety $Z$,  we may
construct another complex algebraic variety $\oZ$
called the (abstract) \emph{complex conjugate} variety.   
If $Z$ is affine, then $\oZ$ is the unique affine variety
such that \[R(\oZ)=\ovl{R(Z)}\]  where 
$\ovl{R(Z)}$ is the $\C$-algebra 
which is complex conjugate to
$R(Z)$; i.e., $\ovl{R(Z)}$ has the same underlying
$\R$-algebra structure but has the complex conjugate
complex vector space structure. For general varieties,
$\oZ$ is defined in the obvious way by gluing together 
complex conjugate affine  opens.

If $f:Z\to Z'$ is a morphism of  complex varieties then
the complex conjugate
map $\of:\oZ\to\ovl{Z'}$ defined by
\[\of(p)=\ovl{f(\ovl{p})}\] is also a morphism.

The construction of $\oZ$ from $Z$ is functorial in the
usual ways and  commutes with products.
We have natural identifications $\ovl{TZ}=T\oZ$
and $\ovl{T^*Z}=T^*\oZ$ for the holomorphic  tangent
and cotangent  bundles. Also pullback of differential
forms and pushforward of vector fields commutes with
taking the complex conjugate.

Consider the natural map
\[Z\to Z\times\oZ,\qquad z\mapsto (z,\oz)\]
 This embeds $Z$ as a real form of  $Z\times\oZ$ with
respect to the real structure defined by
$\ovl{(u,v)}=(\ovl{v},\ovl{u})$.
Thus in particular,   $Z$ itself has a
canonical structure of real algebraic variety.
This amounts to ``forgetting" part of the complex
algebraic variety structure.   Notice that $Z$ and
$\ovl{Z}$ acquire isomorphic real algebraic variety
structures  in this way.

We may write $Z^{real}$ for  $Z$ regarded as real 
variety. If $Z$ is an affine complex variety then
$Z^{real}$ is just the obvious  affine real variety.
Indeed suppose  
$Z\subset\C^n$ is defined by the vanishing of 
$P_1,\cdots,P_m\in\C[z_1,\dots,z_n]$.
We have a natural $\R$-algebra homomorphism
$\C[z_1,\dots,z_n]\to\C[x_1,y_1,\dots,x_n,y_n]$, say
$P\mapsto P'$, defined
by setting   $z'_j=x_j+iy_j$.  
Then $Z^{real}\subset\R^{2n}$ is the closed real 
algebraic subvariety defined by the vanishing of  
the real and imaginary parts
$\Re P'_1,\Im P'_1,\dots,\Re P'_m,\Im P'_m
\in\R[x_1,y_1,\dots,x_n,y_n]$.

If $Z$ has a real structure $\kappa$, then the map
\[Z\to\oZ ,\qquad   z\mapsto\ovl{\kappa(z)}\]
is  an isomorphism of complex
algebraic varieties .

\subsection{Tangent Spaces, Dimension, and 
Smoothness} Let $v$ be a point of an irreducible
algebraic set $V\subset\R^n$. The \emph{ Zariski 
tangent space $T_vV$ at $v$} may be defined as the
linear subspace of $\R^n$ given by
\[T_vV=\{x\in\R^n\,|\, (\hbox{grad} P|_v)\cdot x=
0\hbox{ for all } P\in I(V)\}\]
The  dimension  $d_v=\dim T_vV$ is generically the 
same  over $V$
(i.e., is the same over some  Zariski open dense set
of 
$V$). This common value of $d_v$ is called the 
\emph{dimension}
$d_V$ of $V$. A point $v\in V$ is a \emph{smooth point}
if  $d_v=d_V$. The set $V^{reg}$ of smooth points is 
Zariski open dense in $V$. $V$ is a \emph{smooth variety}
if $V=V^{reg}$.

These notions  pass immediately to affine real 
algebraic varieties and then are purely local.
These notions then pass 
to general real algebraic varieties as the latter are
obtained by gluing of affine opens.  In particular then
the notions of Zariski tangent space and smooth point
are purely local. In the usual way one defines  \emph{\'etale
maps} of  real algebraic varieties.
 
If $Z$ is smooth, then,  in the context of  
\S\ref{ssec:real},  $\Zr$ is a smooth real form of $Z$
(in particular $\Zr$ is irreducible). This follows by
observing that at each point $z\in\Zr$ the complex
Zariski tangent space
$T_zZ$ acquires a real structure and then
the real points form the tangent space to the real 
submanifold  $\Zr$.

A smooth real algebraic variety $X$ has a natural
structure of real analytic manifold,  just as a smooth
complex algebraic variety has a  natural structure of 
complex analytic manifold.  In particular  $X$ has a
larger topology, often called the \emph{strong or
Euclidean  topology}, which refines the Zariski
topology. On $\R^n$, this is just the usual  Euclidean
topology.

Now $X$, while connected in the Zariski topology
(since it is irreducible), may well fail to be connected
in the Euclidean topology. This typically happens when
taking real forms. For example, the familiar real form
of $\C^*$ is $\R^*$. Fortunately,   the
individual Euclidean connected components  have a
natural structure, namely each is a semi-algebraic 
real analytic submanifold. In fact, each component  
is a  \emph{Nash manifold}. We develop this notion in 
the rest of this Appendix.  The starting point is
semi-algebraic sets.

\subsection{Real Semi-Algebraic Sets and Maps}

A subset $S\subset\R^n$
is a  real \emph{semi-algebraic set} if $S$ is a finite 
union of  sets of the form:
\[\{x\in\R^n\,|\, P_1(x)=\cdots=P_m=0 \hbox{ and }
Q_1(x),\dots,Q_m(x)>0\}\]
where $P_i,Q_j\in \R[x_1,\dots,x_n]$.

Suppose $S\subset\R^n$ and 
$T\subset\R^m$ are semi-algebraic sets. 
A map $\phi:S\to T$ is a \emph{semi-algebraic map} if
the graph of $\phi$ is a semi-algebraic
set in $\R^{n+m}$. 
A semi-algebraic map
$f:S\to\R$ is called a \emph{semi-algebraic function}.
It follows that $\phi$ is semi-algebraic if and only if
all the component functions $\phi_1,\dots,\phi_m$ 
are semi-algebraic. 

Notice that a regular map of real algebraic sets
in Euclidean space 
is in particular a semi-algebraic map of 
semi-algebraic sets.

Semi-algebraic sets arise inevitably in the study
of real algebraic sets.
Indeed the image of an algebraic set under a regular
mapping,  even a 
linear projection of Euclidean space, is in general
only semi-algebraic. Also the connected components
(in the Euclidean topology) of an algebraic set
are generally  only semi-algebraic.

The Tarski-Seidenberg Theorem says that under a
semi-algebraic map,
the image of a semi-algebraic set is semi-algebraic.
Another result  says that 
a semi-algebraic set
has finitely many  connected components  (in the 
Euclidean topology) and  each such component  is
semi-algebraic  (see \cite[Th. 2.4.5, pg 31]{BoCR}). 

Next we define semi-algebraic sets in varieties. 
Let $(X,\A_X)$ be a real algebraic variety
and let $S\subset X$.
If $X$ is affine, then we call $S$
\emph{semi-algebraic} if for one (and hence every)
closed embedding $\phi:X\to\R^n$ of real algebraic
varieties, the set $\phi(S)$ is semi-algebraic in $\R^n$.
Now for general $X$ we call $S$ \emph{semi-algebraic}
if for every affine open $U\subset X$ (or equivalently,
for every member $U$ of some affine open cover of
$X$), the set $S\cap U$ is semi-algebraic in $U$.

Notice that if $W\subset X$ is a locally closed
subvariety,
then $S\cap W$ is semi-algebraic in $W$. 

Now, generalizing the definition above, if
$S\subset X$ and $T\subset Y$ are semi-algebraic 
sets in real algebraic varieties, then 
a map $\phi:S\to T$ is \emph{semi-algebraic}
if the graph of $\phi$ is semi-algebraic
in $X\times Y$.
It is routine to check that the Tarski-Seidenberg
theorem is still true in this setting.

A  Euclidean  open semi-algebraic set, and so in 
particular  a  Euclidean  connected component,   of a
smooth real algebraic variety is a  real analytic 
submanifold.

An easy, but important observation is the following:
if $X$ is a real algebraic variety and
$\phi\in A(X)$ is such that $\phi$ takes both positive 
and negative values on $X$ then the set 
\[S=(\phi>0)\subset X\] is semi-algebraic in $X$
(but not algebraic).

\subsection{Nash Functions and Nash Manifolds}
\label{ssec:Nash}
Suppose that $S$ is a (Euclidean) open
semi-algebraic set in a smooth irreducible algebraic 
set $V\subset\R^n$. 

A real analytic function $f:S\to\R$ is called a 
\emph{Nash function} if $f$ satisfies 
the following two equivalent conditions: 

(i) $f$ is
algebraic over the  algebra $P(V)$ of polynomial
functions and 

(ii) $f$ is semi-algebraic. 

The  Nash functions form a Noetherian $\R$-algebra
$N_V(S)$ algebraic over $P(V)$, and furthermore
$N_V(S)$ is integrally closed if $S$ is Euclidean 
connected --- see \cite{BoE}.

From now on assume, more generally,  that  $S$ is a
semi-algebraic real analytic smooth
submanifold    of a smooth  real algebraic variety $X$.

If $X$ is affine, then the  definition  above of Nash 
function on $S$ and the equivalence of the two
conditions go over immediately   as soon
as we  replace  (i) by the condition:  

(i$'$) $f$ is algebraic over $A(X)$. 

\noindent (In the case $X=V$, this is consistent  with the
previous definition as $A(V)$ is algebraic over
$P(V)$) The  Nash functions on $S$ form a Noetherian 
$\R$-algebra $N_X(S)$ which is algebraic over $A(X)$ 
and, if  $S$ is Euclidean connected,  integrally closed. 

Now we can treat the  case where $X$ is
not necessarily  affine. A real analytic function 
$f:S\to\R$ is a \emph{Nash function} if $f$ satisfies 
the following two equivalent conditions: 

(i) for each affine open $U\subset X$ (or equivalently,
for every member $U$ of some affine open cover of
$X$), the restriction
$f|_{S\cap U}$ is algebraic over $A(U)$ and 

(ii) $f$ is semi-algebraic. 

\noindent It follows from the affine case that the  
Nash functions  form an $A(X)$-algebra
$N_X(S)$  which is   integrally closed if $S$ is Euclidean
connected.

Next we define  the \emph{sheaf  $\N_S$ of Nash 
functions on $S$}. We start with the Euclidean
topology on $S$. The collection $\F_S$ of  
semi-algebraic Euclidean open sets in $S$ is a 
basis of this  topology  (e.g, use small open balls).
If $U\in\F_S$, then we  define 
$\N_{S}(U)=N_X(U)$. This data  determines uniquely the 
sheaf  $\N_{S}$ of $\R$-algebras  on $S$.

The pair $(S,\N_{S})$ is then an example of a \emph{Nash 
manifold}. We will not develop a more general theory of
Nash manifolds  here
as these examples are sufficient for  purposes of 
studying  orbits of real algebraic groups, as explained 
in \S\ref{ssec:last} below. 

In particular, smooth real algebraic varieties are Nash 
manifolds and all real algebraic constructions on them
or among them are Nash in the sense discussed below.

Notice that our constructions on $S$ have nice 
functorial properties. For example,
if $X\subset X'$ is a (locally closed) embedding of
smooth real algebraic varieties, then
$X$ and $X'$ define the same Nash manifold structure 
on $S$.

Now suppose  $S'\subset S$ is such that 
$S'$ is a semi-algebraic real analytic smooth
submanifold of $X$. Then $S'$ with its  
sheaf $\N_{S'}$ of Nash functions, is a
\emph{Nash submanifold} of $S$. 
In particular, each  Euclidean connected component of
$S$   is an   open Nash submanifold.

If $(S,\N_{S})$ and $(T,\N_{T})$ are two Nash 
manifolds, then a morphism of the ringed spaces is called
a   \emph{Nash map} or a \emph{Nash morphism}. Thus a map
$\phi:S\to T$ is Nash if and only if
for each Euclidean open set $V\subset T$, 
$\phi\i(V)$ is Euclidean open in $S$  and $\phi$ defines
an algebra homomorphism 
$\phi^*:\N_{T}(V)\to\N_{S}(\phi\i(V))$ by  pullback
of    functions. A Nash map  $\phi$ is a \emph{Nash isomorphism}  if 
$\phi$ is bijective  and $\phi\i$ is Nash.

A Nash map $\phi:S\to T$ is a \emph{Nash embedding}
if $\phi(S)$ is  a Nash submanifold of $T$ and the 
restricted map
$\phi:S\to\phi(S)$ is a Nash isomorphism. 

In the natural way, we define Nash Lie groups,
Nash group actions, etc.

We can define in the obvious way Nash fibrations and 
Nash coverings of Nash manifolds.  We note that local
triviality in the \'etale topology on real algebraic
varieties implies local triviality in the
Euclidean topology.  

Then in particular  we get the notion of a Nash vector bundle 
over a Nash manifold and the space of Nash sections. If $X$ is
a Nash manifold then the tangent and  cotangent bundles of
$X$ have natural Nash bundle structures. Consequently, for
any tensor field $\eta$ on $X$, such as a  vector field, a
differential form, a metric or a complex structure, we define
$\eta$ to be
\emph{Nash} if the corresponding section of the bundle
$TX^{\otimes r}\otimes T^*X^{\otimes s}$ is Nash.
This gives  notions of Nash symplectic manifold, 
Nash Riemannian manifold, Nash complex manifold,
Nash Kaehler manifold, Nash hyperkaehler manifold, etc.

If  $X$ is a totally real Nash submanifold of a smooth complex
algebraic variety $Z$ such that $\dimR X=\dimC Z$, then
we say that $Z$ is a  \emph{Nash complexification} of $X$.
 A stronger condition on $X$ is that $X$ is a Euclidean
connected component of the fixed-point set $Z^\kappa$
for some real structure $\kappa$ on $Z$. Then we say also
that $X$ is a \emph{real form} of $Z$.
This  extends our definition of  real form from 
\S\ref{ssec:real} .

\subsection{Orbits of Real Algebraic Groups}\label{ssec:last}
We consider now real algebraic groups $\GRR$ that arise 
in the following way. Let $G$ be a  Zariski  connected
complex algebraic group defined over $\R$ with group
$\GRR$ of real points. We assume as usual that $G$ is a
complex affine  algebraic variety; then $\GRR$ is a real
affine algebraic variety. For example, 
compact Lie groups  arise  in this way.

Now $\GRR$ is Zariski connected but in
general not Euclidean connected.
For instance if $G=GL(n,\C)$, $n\ge 1$,
then
$\GRR=GL(n,\R)$ has two connected
components, defined by the sign of the determinant.
The Euclidean connected component $\GR$
of $\GRR$ is a semi-algebraic set in $\GRR$. 
If $G$ is semisimple and simply-connected, then 
$\GR=\GRR$.

Suppose $G$ acts morphically on an (irreducible) 
complex algebraic  variety $X$, i.e., $G$ acts on $X$ 
and the action map
$G\times X\to X$ is a morphism of complex algebraic
varieties. If $X$ and the action (i.e., the action 
morphism) are defined over
$\R$ then $\GRR$ acts morphically on the 
(irreducible) set
$\XRR$ of real points. Each Euclidean connected
component $\XR$ of $\XRR$ is a semi-algebraic set
in $\XRR$.

Each $G$-orbit $G\cdot x$ on $X$ is a smooth locally
closed  (irreducible) complex algebraic subvariety of 
$X$. Hence, if $x\in\XRR$, the set of real points
\[(G\cdot x){(\R)}=(G\cdot x)\cap\XRR\]
is a smooth locally closed 
(irreducible) real algebraic subvariety of $\XRR$,
and hence is a finite union of Euclidean connected 
components of the same dimension. These components 
are then semi-algebraic sets and moreover are Nash
submanifolds.

On the other hand, 
by the Tarski-Seidenberg Theorem,
the orbits  $\GRR\cdot x$ and   $\GR\cdot x$ 
are semi-algebraic sets in $\XRR$.  
In particular $\GR\cdot x$ is a 
component of $(G\cdot x){(\R)}$.

Thus Nash manifolds are the natural objects in this 
setting. Finally we give an example of how Nash
isomorphisms can arise.
Consider the standard action of $G=S0(3,\C)$ on 
$\C^3$ as the special orthogonal group of the 
quadratic  form $x^2+y^2-z^2$ where $x,y,z$ are  real
linear coordinates.  The subset
\[X=(x^2+y^2=z^2)-\{(0,0,0)\}\]  is a $G$-orbit.

But $\XRR$ has two Euclidean connected components 
defined by the sign of $z$. Let $\XR$ be the component
where $z>0$.  The projection
\[ p:\XRR\to\C-\{0\},\qquad p(x,y,z)=x+iy\] 
is a $2$-to-$1$ \'etale real algebraic morphism. The
restricted map
$p_\R:\XR\to\C-\{0\}$ defined by $p$
is a Nash isomorphism. Indeed the inverse map is
\[ \C-\{0\}\to\XR,\qquad 
x+iy\mapsto \left(x,y,\sqrt{x^2+y^2}\right)\]


\begin{thebibliography}{99}
 
\bibitem[B1]{B1}  
R. Brylinski, Quantization of the 4-Dimensional
Nilpotent Orbit of $SL(3,\R)$, {\it Can. J. Math.} {\bf 49} (5),
1997, 916-943


\bibitem[B2]{B2}
 R. Brylinski,
Geometric Quantization of Real Minimal
Nilpotent Orbits, preprint 1998


\bibitem[BK1]{bkHam}  R. Brylinski, B. Kostant,  
Nilpotent orbits, normality and Hamiltonian group actions.  
{\it J. Amer. Math. Soc.} {\bf 7} (1994), no. 2, 269--298.

\bibitem[BK2]{bkLagr}
R. Brylinski, B. Kostant,
Lagrangian models of minimal representations of       
$E\sb 6$, $E\sb 7$ and $E\sb 8$,  Functional analysis on the 
eve of the  21st century, 
 Vol. 1 (New Brunswick, NJ, 1993), 13--63, Progr. Math., 131, 
Birkh\"auser Boston, Boston, MA, 1995.

\bibitem[BeR]{BeR}
   R. Benedetti and J-J. Risler, 
   {\it Real Algebraic and Semi-Algebraic Sets},
   Hermann, Paris, 1990.




\bibitem[Bi]{Bi}  O. Biquard, Twisteurs des orbites
coadjointes et m\'etriques
hyper-pseudok\"ahl\'er\-iennes,  Preprint, Ecole
Polytechnique, 1997,  to appear in Bull. Soc. Math. France


\bibitem[BoCR]{BoCR}
  J. Bochnak, M. Coste and M-F. Roy,
  {\it G\'eom\'etrie Alg\'ebrique R\'eelle},
  Ergeb. Math., vol. 12, Springer-Verlag,  1987.


\bibitem[BoE]{BoE}
  J. Bochnak and G. Efroymson,  
  Real Algebraic Geometry and the Hilbert 17th Problem,
  {\it Math. Ann.} {\bf 251} (1980), 213-241.

\bibitem[H]{H}  
     N.J. Hitchin, Metrics on moduli spaces,
{\it Proc. Lefschetz Centennial Conf., Mexico City 1984.
Contemp. Math.}  vol. {\bf 58}, Part I, Amer. Math. Soc.
(1986)
      

\bibitem[HKLR]{HKLR}
     N.J. Hitchin, A Karlhede, U. Lindstrom, and M. Rocek,
     Hyperkaehler metrics and supersymmetry,
     {\it Comm. Math. Phys.} {\bf 108} (1987), 535-559.


   

\bibitem[Kra]{Kra}
   S.G. Krantz,  {\it Function Theory of Several Complex
  Variables}, John Wiley and Sons, 1982


\bibitem[Kr]{Kr2} 
   P. B. Kronheimer, Instantons and the geometry of the nilpotent
   variety, {\it Jour. Diff. Geom.}  {\bf 32} (1990), 473-490.
 

\bibitem[Se]{sek}  
  J. Sekiguchi, Remarks on real
  nilpotent orbits of a symmetric  pair,  {\it J. Math. Soc. 
  Japan}  {\bf 39} (1987), 127-138.

\bibitem[Sw]{Sw} 
  A.  Swann,  Hyper-K\"ahler and quaternionic K\"ahler
  geometry, {\it  Math. Ann.} {\bf 289} (1991), no. 3,
  421-450.

\bibitem[Ve]{Ve}    M. Vergne, 
 Instantons et correspondance de Kostant-Sekiguchi, 
 {\it C.R. Acad. Sci. Paris} {\bf 320} (1995), Serie 1, 901-906.
 
 

\end{thebibliography}
\end{document}